\documentclass[10pt,leqno]{article} 
\usepackage{amsmath,amsfonts}

\newsymbol \flechas 1313

\newcommand{\bl}{\begin{lemma}}
\newcommand{\bp}{\begin{proposition}}
\newcommand{\bt}{\begin{theorem}}
\newcommand{\bc}{\begin{corollary}}

\newtheorem{proposition}{Proposition}[section] 
\newtheorem{lemma}[proposition]{Lemma} 
\newtheorem{corollary}[proposition]{Corollary} 
\newtheorem{theorem}[proposition]{Theorem} 

\newcommand{\el}{\end{lemma}}
\newcommand{\ep}{\end{proposition}}
\newcommand{\et}{\end{theorem}}
\newcommand{\ec}{\end{corollary}}

\newcommand{\bL}{\begin{Lemma}}
\newcommand{\bP}{\begin{Proposition}}
\newcommand{\bT}{\begin{Theorem}}
\newcommand{\bC}{\begin{Corollary}}

\newtheorem{Proposition}{Proposition}[proposition] 
\newtheorem{Lemma}[Proposition]{Lemma} 
\newtheorem{Corollary}[Proposition]{Corollary} 
\newtheorem{Theorem}[Proposition]{Theorem}

\newcommand{\eL}{\end{Lemma}}
\newcommand{\eP}{\end{Proposition}}
\newcommand{\eT}{\end{Theorem}}
\newcommand{\eC}{\end{Corollary}}

\newcommand{\prg}{ \setcounter{Proposition}{0}                     
                     \addtocounter{proposition}{1}   
                        \paragraph{\theproposition} }

\newcommand{\prgg}{  \addtocounter{Proposition}{1}   
                      \paragraph{\theProposition } }

\newcommand{\TO}{\longrightarrow}
\newcommand{\nt}{\noindent}  
\newcommand{\wt}{\widetilde} 
\newcommand{\refp}[1]{(\ref{#1})}

\newcommand{\pro}{\noindent {\it Proof. \/}}
\newcommand{\qed}{\hfill $\clubsuit$ \medskip}
\newcommand{\len}{\mathop{\rm len \, } \nolimits}
\newcommand{\Ker}{\mathop{\rm Ker \, } \nolimits}
\newcommand{\Coker}{\mathop{\rm Coker \, } \nolimits}
\newcommand{\Ima}{\mathop{\rm Im \, } \nolimits}
\newcommand{\inte}{\mathop{\rm int \, } \nolimits}
\newcommand{\codim}{\mathop{\rm codim \, } \nolimits}
\newcommand{\dimc}{\mathop{\rm dimc \, } \nolimits}
\newcommand{\id}{\mathop{\rm id \, } \nolimits}
\newcommand{\proj}{\mathop{\rm projection \, } \nolimits}
\newcommand{\rond}{\hbox{\small\hbox{$\circ$}}}
\newcommand{\Tor}{\mathop{\rm Tor \, } \nolimits}
\newcommand{\Id}{\mathop{\rm id \, } \nolimits}

\newcommand\chii{\raise2pt\hbox{$\chi$}}      
\newcommand\phii{{\raise2pt\hbox{$\varphi$}}}  
\newcommand\Om{\Omega}   
\newcommand\om{\omega}  

\newcommand\C{\mathbb C}
\newcommand\N{\mathbb N}
\newcommand\R{\mathbb R}

\renewcommand\S{\mathbb S}
\newcommand\T{\mathbb T}
\renewcommand\P{\mathbb P}

\newcommand\Z{\mathbb Z}

\newcommand{\sbat}{{\S}^{^1}}

\newcommand{\bi}[2]{{#1}^{^{#2}}}
\newcommand{\hiru}[3]{{#1}^{^{#2}}{\left( #3 \right)}}
\newcommand{\hirut}[3]{{\wt{#1}}^{^{#2}}{\left( #3 \right)}}
\newcommand{\hiruv}[3]{{#1}^{^{#2}}_{_v}{\left( #3 \right)}}
\newcommand{\hirue}[3]{{#1}^{^{#2}}_{\sbat}{\left( #3 \right)}}

\newcommand{\dos}[2]{{\mathcal #1}(#2)}

\catcode`\Š=\active\defŠ{\"a}        % option u, then  a
\catcode`\'=\active\def'{\"e}        % option u, then  e
\catcode`\•=\active\def•{\"{\i}}     % option u, then  i
\catcode`\š=\active\defš{\"o}        % option u, then  o
\catcode`\Ÿ=\active\defŸ{\"u}        % option u, then  u
\catcode`\Ø=\active\defØ{\"y}        % option u, then  y
\catcode`\€=\active\def€{\"A}        % option u, then  A
\catcode`\…=\active\def…{\"O}        % option u, then  O
\catcode`\†=\active\def†{\"U}        % option u, then  U
\catcode`\‡=\active\def‡{\'a}        % option e, then  a
\catcode`\Ž=\active\defŽ{\'e}        % option e, then  e
\catcode`\'=\active\def'{\'{\i}}     % option e, then  i
\catcode`\—=\active\def—{\'o}        % option e, then  o
\catcode`\œ=\active\defœ{\'u}        % option e, then  u
\catcode`\ƒ=\active\defƒ{\'E}        % option e, then  E
\catcode`\ˆ=\active\defˆ{\`a}        % option `, then  a
\catcode`\=\active\def{\`e}        % option `, then  e
\catcode`\"=\active\def"{\`{\i}}     % option `, then  i
\catcode`\˜=\active\def˜{\`o}        % option `, then  o
\catcode`\=\active\def{\`u}        % option `, then  u
\catcode`\Ë=\active\defË{\`A}        % option `, then  A
\catcode`\‹=\active\def‹{\~a}        % option n, then  a
\catcode`\–=\active\def–{\~n}        % option n, then  n
\catcode`\›=\active\def›{\~o}        % option n, then  o
\catcode`\Ì=\active\defÌ{\~A}        % option n, then  A
\catcode`\"=\active\def"{\~N}        % option n, then  N
\catcode`\Í=\active\defÍ{\~O}        % option n, then  O
\catcode`\‰=\active\def‰{\^a}        % option i, then  a
\catcode`\=\active\def{\^e}        % option i, then  e
\catcode`\"=\active\def"{\^{\i}}     % option i, then  i
\catcode`\™=\active\def™{\^o}        % option i, then  o
\catcode`\ž=\active\defž{\^u}        % option i, then  u   
                                
\title{\Large\bf MINIMAL MODELS \\ FOR \\ NON-FREE CIRCLE ACTIONS\footnote{Accepted for publication in Illinois Journal of Mathematics}}

\author{ Agust' Roig\thanks{Departament de Matemˆtica Aplicada.
E. T. S. E. I. B.  Universitat Politnica de Catalunya. Diagonal 647. 08028 Barcelona - Spain. 
roig@ma1.upc.es}  \and   Martintxo Saralegi-Aranguren\thanks{Laboratoire GŽomŽtrie-Algbre. 
 FacultŽ Jean Perrin. UniversitŽ d'Artois.  Rue Jean Souvraz SP 18.  62307 Lens Cedex - France.  
saralegi@euler.univ-artois.fr}}

\begin{document}  

\maketitle

\begin{abstract} 
Let $\Phi \colon \sbat \times M \to M$
be a smooth action of the unit circle $ \sbat$ on a
manifold $M$. In this work, 
we compute the minimal model of $M$ in
terms of the orbit space $B$ and the fixed
point set $F\subset B$, as a dg-module over the
Sullivan's minimal model of $B$.
\end{abstract}

\thispagestyle{empty}

The question we treat in this work is the following:
given a smooth action $\Phi \colon \sbat \times M \to
M$, is it possible to construct a model of $M$ using
just basic data?
 The answer is well known when the fixed
point set is empty (in particular, when the action is
free). A dgc algebra model of
$M$ is  given by a Hirsch extension of the dgc algebra
Sullivan minimal model $\dos{A}{B}$ of the orbit space
$B$
\begin{equation}
\label{galdera}   
\dos{A}{B} \otimes \Lambda(x),
\end{equation}  where
the degree of $x$ is 1  and $dx$ defines the Euler
form of the action (see for example
\cite{GMo}). This formula does not apply when the fixed
point set $F$ is not empty. Roughly speaking this
happens because the Euler form does not live on
$B$.

Our answer to the above question is a minimal model
of the deRham dgc algebra $\hiru{\Om}{}{M}$ of
$M$, which is a dg module over the Sullivan
minimal model $\dos{A}{B}$ of $B$.  Such structure is
associated to
$M$ by means of the canonical projection
$\pi \colon M \to B$. We prove that the minimal
model of $M$, as an $\dos{A}{B}$-dg module,  is the
graded cone
\begin{equation}
\label{model}
\dos{M}{M} = \dos{A}{B} \oplus_{e'} \dos{M}{B,F} ,
\end{equation}   
where $\dos{M}{B,F}$ is a sort of relative
minimal model of the pair $(B,F)$ and $ e'\colon
\dos{M}{B,F}
\to
\dos{A}{B}$ is a degree 2 map. This map is
determined by the Euler class of the action and it
will be described below. We also prove that, for $F =\emptyset$, the 
formul¾
\refp{galdera} and \refp{model} coincide.

There are some algebraic invariants of $M$ and $F$
that are closely related: PoincarŽ characteristic,
localization, rational homotopy, \ldots We add another
item  to this list: the minimal model of $M$ and $F$.
In fact, considering the $\dos{A}{B}$-dg module
structure associated to $F$ by means of the natural
inclusion
$\iota \colon F \hookrightarrow B$, we prove that the
minimal model of $F$ as an $\dos{A}{B}$-dg module is
the graded cone
$$
\dos{M}{F} = \dos{A}{B} \oplus_{i'} \dos{M}{B,F} ,
$$  
where $ i'\colon \dos{M}{B,F}
\to
\dos{A}{B}$ is a degree 0 map. This map is
determined by $\iota$ and will be described
below. Observe that the minimal models $\dos{M}{M}$
and $\dos{M}{F}$ have, as
$\dos{A}{B}$-graded modules, the same basis except
for a shift by 2.

We are also interested in the Borel space $M_{\sbat} =
M \times_{\sbat} \S^\infty$, where an
$\dos{A}{B}$-dg module structure is defined  by
means of the canonical projection $p \colon M_{\sbat}
 \to B$. We prove that the minimal model of
$M_{\sbat}$ is  the graded cone
$$
\dos{M}{M_{\sbat}} = \left[ \dos{A}{B} \otimes
\Lambda(e) \right]
 \oplus_{q'}
\left[\dos{M}{B,F}\otimes \Lambda(e) \right],
$$ 
where $\deg e = 2$, $de=0$ and $q'(b\otimes e^n) =e'(b)
\otimes e^n + i'(b)
\otimes e^{n+1}$. This formula implies that the
equivariant cohomology $\hirue{H}{i}{M}$ (i.e.
the cohomology of
$M_{\sbat}$) can be computed using just basic
data by means of the long exact sequence
 $$
\cdots \to \left[ \hiru{H}{}{B} \otimes \Lambda (e)
\right]^{^i} \to
\hirue{H}{i}{M} \to \left[ \hiru{H}{}{B,F} \otimes
\Lambda (e) \right]^{^{i-1}}
 \stackrel{(q')^*}{\TO} \left[ \hiru{H}{}{B} \otimes
\Lambda (e) \right]^{^{i+1}}
\to \cdots .
$$
Moreover,
when the Euler class vanishes, we prove that the
equivariant cohomology of
$M$ is just 	\newline
$
\hiru{H}{}{B} \oplus \left[ \hiru{H}{}{F}
\otimes \Lambda^+ (e)\right] .
$
We also translate some classic results
(Localization Theorem, equivariant
formality, \ldots) in terms of basic data.
In connection with the minimal models of $M$,  $M_{\sbat}$ and $F$ 
the 
reader can consult \cite{A} and \cite{AP}.

Let us illustrate these results with the suspension of
the Hopf action on  $\S^3$. The north and south poles
of the total space $\S^4$ are the fixed points of the
action and the orbit space is $\S^3$. We work in the
category of
$\Lambda (a)$-dg modules, where $\deg a = 3$. From the
above formul¾ we get:
$$
\begin{array}{rcrcl}
\dos{M}{\S^4}&=&\Lambda (a)& \otimes & \R 
\{ 1, c_n \ / \  n\in \N \} ,
\hbox{ with} \\[,2cm]  &&&& \deg c_n=2[\frac{n+3}{2}],
d c_0 =a, dc_1=0,  dc_{n+2}=a\cdot c_n.\\[,5cm]
\dos{M}{\S^0}&=&\Lambda (a)& \otimes & \R 
\{ 1, \gamma_n \ / \  n\in \N \} ,
\hbox{ with} \\[,2cm] & &&& \deg
\gamma_n=2[\frac{n+1}{2}], d
\gamma_0 =0, d\gamma_1=a, d\gamma_{n+2}=a\cdot
\gamma_n.\\[,5cm]
\dos{M}{\S^4\otimes_{\sbat} \S^\infty}&=&\Lambda (e,a)
&\otimes & \R 
\{ 1, c_n \ / \  n\in \N \} ,
\hbox{ with } \\[,2cm] &&&&  \deg
c_n=2[\frac{n+3}{2}],  d c_0 =a, dc_1=e \cdot
a, dc_{n+2}=a\cdot c_n. 
\end{array}
$$ 
The minimal  model $\dos{M}{\S^4}$ (resp.
$\dos{M}{\S^0}$, resp. $\dos{M}{\S^4\otimes_{\sbat}
\S^\infty}$) is a free $\Lambda (a)$-graded module over
the cohomology $\hiru{H}{*}{Y_\pi}$ of the homotopy
fiber $Y_\pi$ of $\pi$  (resp. $Y_\iota$, resp.
$Y_p$). So, we find the following relations between
the PoincarŽ polynomials of these spaces:
$$ P_{Y_\pi} = 1 - t^2 + t^2P_{Y_\iota}  = (1 - t^2) 
P_{Y_p}.
$$ We prove that these relations are generic if $B$ is
simply connected and of finite type.

The main geometric tool used in this work are
Verona's controlled forms
\cite{V}.  In fact, when the set of fixed points $F$
is not empty the orbit space
$B$  is not a regular manifold but a singular one,
more precisely a stratified pseudomanifold. For a such
space $Z$, Verona proved that the complex of
controlled forms
$\hiruv{\Om}{*}{Z}$ compute the cohomology of $Z$. We
prove  more: that the minimal models $\dos{A}{Z}$ and
$\dos{M}{Z}$ can be computed using controlled forms.
It is important to notice that the Euler form is not a
controlled form, nevertheless  it appears in this
context as a morphism of
$\dos{A}{B}$-dg modules 
$e \colon \hiruv{\Om}{*}{B,F} \to
\hiruv{\Om}{*+2}{B}$ (cf. \cite{HS}). In the writing
of
$\dos{M}{M}$ (resp. $\dos{M}{F}$) the operator $e'$
(resp. $i'$) is a model of $e$ (resp. of the inclusion
$i \colon
\hiruv{\Om}{*}{B,F} \hookrightarrow
\hiruv{\Om}{*}{B}$).

The starting point of the work is the observation that
the cohomology of $M$ can be computed by the graded
cone
$
\hiruv{\Om}{*}{B} \oplus_e
\hiruv{\Om}{*}{B,F}.
$ This formula also applies to semi-free actions of
$\S^3$ \cite{S2}. So, all the results of this work
extend to this kind of actions. A similar formula
appears when one deals with an isometric action
$\Phi \colon \R \times M \to M$,  considering on $B$
controlled basic forms instead of controlled forms 
 \cite{S2}. Again, we conclude that the results of
this work apply to isometric flows. In particular, we
get the inequality
$$
\hiru{H}{r-1}{(M,F)/{\cal F}} +\sum_{i=0}^{\infty}
\dim \hiru{H}{r+2i}{F} \leq 
\sum_{i=0}^{\infty} \dim \hiru{H}{r+2i}{M},
$$ when the flow is not trivial.

On the algebraic side, we develop in some extend the
Theory of dg minimal modules. This kind of minimal
objects was previously studied by the first author
(cf. \cite{R1}, \cite{R2}, \cite{R3}) and
independently by Kriz and May (cf. \cite{KM}).

 \bigskip

The organization of the work is as follows. First
section is devoted to present the algebraic tools we
need to work with $A$-dg modules. In the second section
we present the singular spaces we find when we deal
with circle actions. Controlled forms are introduced
in the third section. The main result of this paper is
proved in the forth and last section. Four technical
Lemma are proven in the Appendix.

A manifold is considered to be connected, without
boundary and smooth (of class $C^\infty$), unless
otherwise is stated. The field of coefficients is $\R$.

We thanks Yves Felix, Steve Halperin, Pascal
Lambrechts, Vicente Navarro Aznar and Daniel TanrŽ for
their useful comments. 
The authors thank the referee for their useful indications.

\bigskip

\nt {\bf Convention}. Through all this paper we will
adhere to the following convention: minimal models in
the category of dgc algebras will be denoted by
$\dos{A}{-}$. Minimal models in the category of dg
modules will be denoted by $\dos{M}{-}$.

\section{Dg-Module minimal models}

In this section we will develop the algebraic machinery  necessary to prove
Theorem
\ref{nagusia}: we define what is meant by a {\it minimal factorization of a
morphism} of
$A$-dg modules, prove its existence and uniqueness and a result concerning maps
induced between them.

\prg {\bf \mbox{\boldmath $A$}-dg modules.} Let  $A$ be a dgc algebra. An
$A$-{\em dg module} $M$ is a graded   vector space together with a product $A
\otimes M
\longrightarrow M $ and a differential $d\colon M\longrightarrow M$ of degree
$+1$ which satisfies Leibniz rule. Both  graduations of $A$ and
$M$ are over the non negative integers. A {\em quasi-isomorphism} ({\em quis}) is
an
$A$-dg module morphism which induces an isomorphism in cohomology.

Let us begin with an immediate generalization of the cone of a morphism of
complexes in the category of $A$-dg modules. Let
$\phii \colon M \to N$ be a morphism of $A$-dg modules of degree $p\in \Z$. This is
the same as a degree 0 morphism of
$A$-dg modules
$$
\phii \colon M[-p] \to N
$$ where $M[-p]$ means the $A$-dg module $M$ shifted by $-p$; it is graded by
$M[-p]^n=M^{n-p}$ and the product $\mu_{M[-p]}
\colon A \otimes M[-p] \to M[-p]$ and the differential
$d_{M[-p]} \colon M[-p] \to M[-p]$ change signs according to
$$
\begin{array}{rcl}
\mu_{M[-p]} (a\otimes m) &=& (-1)^{|a|\cdot p} \mu_{M} (a\otimes m)\\ d_{M[-p]} (
m) &=& (-1)^p d_M ( m).
\end{array}
$$ We shall denote by $N\oplus_\phii M$ the $A$-dg module graded by
$$ N\oplus_\phii M = N \oplus M[1-p]
$$ and with product and differential given by the formulas:

$$
\begin{array}{rcl} a \cdot 
\left(\begin{array}{c}  y \\ x
\end{array}\right)  & = &
\left(\begin{array}{c}  a \cdot y \\ (-1)^{^{|a|\cdot (p-1)} }a \cdot x
\end{array}\right)
\\[1cm]
\left(
\begin{array}{cc} d & \phii\\ 0& (-1)^{p-1} d
\end{array}
\right) 
\left(\begin{array}{c}  y \\ x
\end{array}
\right) &=&
\left(\begin{array}{c}  d y + \phii x \\ (-1)^{^{p-1}} d x
\end{array}\right),
\end{array}
$$ where $a \cdot x$ denotes $\mu_M (a \otimes x)$. Finally, when we say that the
sequence of $A$-dg module morphisms
$$ 0 \TO M \stackrel{\varphi}{\TO} N \stackrel{\psi}{\TO} P \TO 0
$$ is {\em exact}, we simply mean that
$ 0 \TO M^n \stackrel{\varphi^n}{\TO} N^n \stackrel{\psi^n}{\TO} P^n \TO 0
$ are usual exact sequences of $A^0$-modules for all $n$.

\prgg {\bf Remark}. A short exact sequence of $A$-dg modules as above is the
same as the {\em quis} of $A$-dg modules
$$ N \oplus_\phii M \stackrel{(\psi,0)}{\TO} P.
$$

\bL
\label{112} Given a morphism $\phii \colon M \to N$ of $A$-dg modules of degree
$p$,  we have a short exact sequence of $A$-dg modules
$$ 0 \TO N \stackrel{\left( ^1_0 \right)}{\TO} N
\oplus_\phii M
\stackrel{(0,1)}{\TO} M[1-p]
\TO 0.
$$
\eL

\pro Obvious. \qed

Associated to the above exact sequence, there is the {\em long exact cohomology
exact sequence} of $\phii$:
$$
\cdots \to \hiru{H}{n}{N}  \stackrel{\left( ^1_0 \right)_*}{\TO} \hiru{H}{n}{N
\oplus_\phii M}
\stackrel{(0,1)_*}{\TO}
\hiru{H}{n+1-p}{M}
\stackrel{\phii_*}{\TO}
\hiru{H}{n+1}{N} \to \cdots.
$$ To end this elementary differential homological algebra, let us point out that a
({\em degree $p$}) {\em homotopy} between two $A$-dg module morphisms $\phii,
\psi \colon M \to N$ of degree $p$,  is an $A$-dg module morphism  $h \colon M \to
N$ of degree
$p-1$, such that
$$ (-1)^p dh + hd = \psi -\phii.
$$ One can verify that this notion of homotopy coincides with the one defined in
\cite{R3} using a path object.

\prg {\bf Minimal \mbox{\boldmath $A$}-dg modules.} (cf.
\cite{Hal}, \cite{Na}, \cite{R1}) Let  $M$ be a $A$-dg module and $n$ a non negative
integer. A {\em degree
$n$ Hirsch extension} of $M$ is an inclusion  of $A$-dg modules
$ M
\hookrightarrow M\oplus (A \otimes V^{^n}) $ in which: 
\begin{itemize} \item[1.] $V$ is a homogeneous vector space of degree
$n$,  
\item[2.] $A \otimes V^{^n}$ is the free $A$-graded module over
$V$, and  
\item[3.]  The differential of $ M\oplus (A \otimes V^{^n})$ is induced
by the differentials of $M$ and $A$ and the choice of a linear map $d\colon
V^{^n} \to M^{^{n+1}}$.
\end{itemize}

A morphism of $A$-dg modules $ M \oplus (A \otimes V^{^n})
\longrightarrow N $ is given by a morphism of $A$-dg modules 
$\phii \colon M\longrightarrow  N$ and  a linear map
$f\colon V^{^n}\longrightarrow N^{^n}$  subjected to the condition $ \phii
\rond d  = d \rond f. $

A {\em minimal KS-extension} of $M$ is an  inclusion of $A$-dg modules $ \iota
\colon M \longrightarrow N $ together with an exhaustive filtration  $\{ N(n,q)
\}_{(n,q)
\in I}$ of $N$, indexed by  $I = \{ (n,q) \in \N \times \N 
\}$ with lexicographical order, such that: 
\begin{itemize}
\item[1.] $N(0,0) = M$,
\item[2.] for
$q > 0$, $N(n,q)$ is a degree $n$ Hirsch extension of
$N(n,q-1)$, and
\item[3.] $N(n+1,0) =  {\displaystyle 
\lim_{\stackrel{\longrightarrow}{q}}} N(n,q)$.  \end{itemize}

\prgg {\bf Remark.} $N$ is therefore of the form $M \oplus (A
\otimes V)$, where $V$ is a bigraded vector space. This kind of objects plays the
r™le of the KS-extensions of \cite{Hal}, noted $B \otimes \Lambda V$, with the
tensor product replaced by the direct sum and the free dg-algebra over $V$
replaced by the free $A$ -dg module over $V$.

\bigskip

 A {\em minimal KS-factorization} of an $A$-dg module morphism
$\phii : M\longrightarrow N$ is a commutative  diagram of
$A$-dg module morphisms

\bigskip

$$ \begin{picture}(100,40)(0,0) 
\put(50,40){\makebox(0,0){$M$}}
\put(0,0){\makebox(0,0){$M \oplus (A \otimes V)$}} 
\put(100,0){\makebox(0,0){$N$}}

\put(38,0){\vector(1,0){50}}
 \put(40,35){\vector(-1,-1){25}} 
\put(60,35){\vector(1,-1){25}}

\put(63,7){\makebox(0,0){$\rho$}}
\put(25,30){\makebox(0,0){$\iota$}}
\put(75,30){\makebox(0,0){$\phii$}}

\addtocounter{equation}{1}
\put(-180,20){\makebox(0,0){(\theequation)}}
\end{picture}                                            $$

\bigskip

\nt in which $\rho $ is a {\em quis} and $\iota $ is a minimal  KS-extension. If $M$ is
the zero $A$-dg module, we speak of {\em minimal KS-modules} and {\em
minimal KS-models}.

 \prg {\bf Models of $A$-dg modules}. Let now (\theequation) be any commutative
diagram of $A$-dg module morphisms. In this situation, we will say that  $\rho $ is
an
$M$-{\em morphism}. If
$\rho $ is also a {\em quis}, we will simply say that it is a
$M$-{\em quis}. A homotopy between two $M$-morphisms which restricted to $M$
is the identity will be called also a
$M$-{\em homotopy}.

\bigskip

Alternatively, we could have said that $\rho $ is a morphism of $ M \backslash {\bf
DGM}(A)$,  {\em the category of $A$-dg modules under}
$M$. So, a minimal KS-factorization is, simply, a minimal model in $ M\backslash {\bf
DGM}(A)$ (see \cite{R3} for the precise statement of this).

\bT  
\label{THI} Let $A$ be a dgc algebra and let $\phii \colon M
\longrightarrow N$ be an $A$-dg module morphism such that
$\phii^{^0}_*
\colon \hiru{H}{0}{M}
 \longrightarrow \hiru{H}{0}{N}$ is a monomorphism. Then there exists a minimal
KS-factorization of $\phii $.  
\eT

\pro See Appendix. \qed

\bC Let $A$ be a dgc algebra and let $N$ be an $A$-dg module. Then there exists a
minimal KS-model of $N$.
\eC

\bT \label{THII}(cf. \cite{H-T}, \cite{R1})  Let $A$ be a dgc algebra and 
$$
\begin{picture}(100,40)(0,0) 
\put(50,40){\makebox(0,0){$M$}}
\put(0,0){\makebox(0,0){$X$}} 
\put(100,0){\makebox(0,0){$M \oplus (A \otimes V)$}}

\put(12,0){\vector(1,0){50}}
 \put(40,35){\vector(-1,-1){25}}
 \put(60,35){\vector(1,-1){25}}

\put(37,7){\makebox(0,0){$\phii$}}
\put(25,30){\makebox(0,0){$\psi$}}
\put(75,30){\makebox(0,0){$\iota$}}

\end{picture}                                            
$$  a commutative diagram of
$A$-dg module morphisms in which  $\iota $ is a minimal KS-extension and
$\phii $ is a {\rm quis}. Then there exists an
$A$-dg module morphism $\sigma \colon M \oplus (A \otimes V)
\longrightarrow X$ such that
$\sigma \iota = \psi  $ and $\phii \sigma = \
\hbox{id}$. 
\eT

\pro See Appendix. \qed

\smallskip

In other words, every $M$-{\em quis} whose target is a minimal KS-extension of
$M$ has a  section which is  also  an
$M$-morphism. This implies the  uniqueness up to
isomorphism of minimal models well known in other categories.

\bC   Two minimal KS-factorizations of the same $M$-morphism are
$M$-isomorphic and the isomorphism is unique up to
$M$-homotopies. \eC

\pro It follows easily from Theorem \ref{THII}, taking into account that
$M\backslash {\bf DGM}(A)$ is a closed model category in which all objects are
fibrant (cf.
\cite{Qui}, \cite[Corollary 2 to Proposition 1.15]{R3}). 
\qed

\smallskip

Particularly, if we take $M=0$, the zero $A$-dg module, we obtain

\bC
\label{uno} Two minimal KS-models of the same $A$-dg module are isomorphic and
the isomorphism is unique up to homotopies.
\eC

Given a morphism $\phii \colon M \to N$  of $A$-dg modules we will need to
construct a model of $N \oplus_\phii M$. This will be done by means of the
following results:

\bC
\label{loose}
 Let $\phii \colon M\longrightarrow N$ be a morphism of
$A$-dg modules and let
$\rho_M
 \colon M'\longrightarrow M$ and $\rho_N \colon N'\longrightarrow N$ be two
minimal models. Then there exists a morphism $\phii'
\colon M'\longrightarrow N'$, unique up to homotopies, that renders commutative
up to homotopy the diagram 

\bigskip

 $$ 
\begin{picture}(50,50)(0,0) 
\put(0,50){\makebox(0,0){$M'$}}
\put(50,50){\makebox(0,0){$M$}} 
\put(0,0){\makebox(0,0){$N'$}}
\put(50,0){\makebox(0,0){$N$}}

\put(0,40){\vector(0,-1){30}} \put(50,40){\vector(0,-1){30}}
\put(12,0){\vector(1,0){30}} \put(12,50){\vector(1,0){30}}

\put(10,25){\makebox(0,0){$\phii'$}}
\put(60,25){\makebox(0,0){$\phii$}}
\put(27,7){\makebox(0,0){$\rho_N$}}
\put(27,57){\makebox(0,0){$\rho_M$}}

\end{picture}                                            
$$  We will loosely say that $\phii'$ is a {\em model} of $\phii$.
\eC

\pro It follows from Theorem \ref{THII} and the model category structure
 \cite[Proposition 1.16]{R2}.  \qed

\bP
\label{suma} Consider the above diagram.   Let $h \colon M' \to N$ be a homotopy
between
$\phii \rond \rho_M$ and $\rho_N \rond \phii'$. Then
$$
\Phi = 
\left(
\begin{array}{cc}
\rho_N & h\\ 0 & \rho_M
\end{array}
\right) 
\colon N' \oplus_{\phii'} M' \TO N \oplus_\phii M
$$ is a {\rm quis}.
\eP

\pro Let us verify that $\Phi$ commutes with differentials.

$$
\begin{array}{rc}
\left(
\begin{array}{cc} d & \phii\\ 0 &(-1)^{^{p-1}} d
\end{array}
\right) 
\left(
\begin{array}{cc}
\rho_N & h\\ 0 & \rho_M
\end{array}
\right)  -
\left(
\begin{array}{cc}
\rho_N & h\\ 0 & \rho_M
\end{array}
\right) 
\left(
\begin{array}{cc} d & \phii\\ 0 &(-1)^{^{p-1}} d
\end{array}
\right)  & = \\[,5cm]
\left(
\begin{array}{cc} d\rho_N -\rho_N' &dh +(-1)^{^{p}} hd - (\rho_N \phii' - \phii
\rho_M) \\[,3cm] 0 &(-1)^{^{p-1}} d \rho_M - (-1)^{^{p-1}}\rho_Md' 
\end{array}
\right). 
\end{array}
$$ And this is zero because $\rho_M$ and $\rho_N$ commute with differentials and
$h$ is a degree $p$ homotopy between $\phii \rond \rho_M$ and $\rho_N \rond
\phii'$. Next, we put $\Phi$ in the following obviously commutative diagram with
exact rows:

$$ 
\begin{picture}(225,70)(0,-10) 

\put(0,50){\makebox(0,0){$0$}}
\put(50,50){\makebox(0,0){$N$}}
\put(100,50){\makebox(0,0){$N \oplus_\phii M$}}
\put(175,50){\makebox(0,0){$M[1-p]$}}
\put(225,50){\makebox(0,0){$0$}}

\put(0,0){\makebox(0,0){$0$}}
\put(50,0){\makebox(0,0){$N'$}}
\put(100,0){\makebox(0,0){$N' \oplus_{\phii' } M'$}}
\put(175,0){\makebox(0,0){$M'[1-p]$}}
\put(225,0){\makebox(0,0){$0$}}

\put(50,10){\vector(0,1){30}} 
\put(100,10){\vector(0,1){30}}
\put(175,10){\vector(0,1){30}} 

\put(10,0){\vector(1,0){30}}
\put(60,0){\vector(1,0){10}}
\put(130,0){\vector(1,0){15}}
\put(205,0){\vector(1,0){10}}

\put(10,50){\vector(1,0){30}}
\put(60,50){\vector(1,0){10}}
\put(130,50){\vector(1,0){15}}
\put(205,50){\vector(1,0){10}}

\put(63,25){\makebox(0,0){$\rho_N$}}
\put(110,25){\makebox(0,0){$\Phi$}}
\put(188,25){\makebox(0,0){$\rho_M$}}

\end{picture}                                            
$$  And so $\Phi$ is a {\em quis} by the long exact sequence of $\phii$ and $\phii'$
and the Five Lemma.
\qed

Finally, we will need the following result concerning minimal dg modules and
graded cones. 

\bL
\label{minimal} Let $\phii \colon M \to A$ be a degree $p$ morphism of $A$-dg
modules, with
$M$ a minimal $A$-dg module.  If $M^{<(1-p)} =0$ then $A \oplus_\phii M$ is  a
minimal $A$-dg module.
\eL

\pro

Let $M= A \otimes V$. Then, as a graded module, $A \oplus M = A \otimes (\R
\oplus V)$. So we can define an exhaustive filtration in $A \oplus M$ as follows. Let
$W(n,q)$ be the $\R$-vector spaces
$$ W(n,q) = \left\{
\begin{array}{ll} 0 & \hbox{ if } n=q=0\\
\R & \hbox{ if } n=0 \hbox{ and } q=1\\ V(n+1-p,q-1) & \hbox{ if } n=0 \hbox{ and }
q>1\\ V(n+1-p,q) & \hbox{ if }  n\not=0
\end{array}
\right.
$$ and put
$$ (A\oplus_\phii M)(n,q) = A \otimes \left(\bigoplus_{(m,r) \leq
(n,q)}W(m,r)\right).
$$ Then, all the inclusions $(A\oplus_\phii M)(n,q-1) \hookrightarrow (A\oplus_\phii
M)(n,q)$ are degree $n$ Hirsch extensions.
\qed

\prgg {\bf Links with Topology}. Before studying the relative case, let us show one
example where these {\em minimal dg-modules} appear in topology. Let $p
\colon E \to B$ be a continuous map between two topological spaces. We have the
induced morphism $p^* \colon A_\R (B) \to A_\R (E)$ between the real algebras of
polynomial forms making $A_\R (E)$ an $A_\R (B)$-dg module. Now, assume that
$B$ is connected, of finite type, with $\pi_1 (B)$ acting trivially on
$\hiru{H}{*}{Y_p} = \hiru{H}{*}{Y_p,\R}$, where $Y_p$ denotes the {\em homotopy
fiber} of
$p$. Then, by the second Theorem of Eilenberg-Moore (see \cite{GM} and
\cite{BG}), we have that 
$$
\hiru{H}{*}{Y_p} \cong \Tor_{A_\R (B)} (\R , A_\R (E)).
$$ By \cite{R3}, this differential torsion product can be computed with a minimal
model of $A_\R (E)$ as an $A_\R (B)$-dg module. Let 
$\dos{M}{E} = A_\R (B) \otimes_\R V$ be this minimal model. Then
$$
\hiru{H}{*}{Y_p} \cong \hiru{H}{*}{\R \otimes_{A_\R (B)} \dos{M}{E}} \cong 
\hiru{H}{*}{\R \otimes_{A_\R (B)} (A_\R (B) \otimes V)} = V
$$ because $\R \otimes_{A_\R (B)} (A_\R (B) \otimes V) = V$ has zero differential
due to minimality. So, the known   Hirsch-Brown model of $E$
$$
\dos{M}{E} =A_\R (B) \otimes \hiru{H}{*}{Y_p} 
$$ is a {\em minimal model} as $A_\R (B)$-dg modules. In particular, if we take $B$
to be a point, we find that $\hiru{H}{*}{E}$ is the minimal model of $A_\R (E)$ as a
$\R$-dg module.

 Extensive use of the minimal Hirsch-Brown model for 
the Borel construction is made in  \cite{AP}.

\prg {\bf Models of couples, existence and uniqueness}. In fact, we will need
something more than simply dg-minimal models over a fixed
$\R$-dg algebra. The process we are going to perform is the following: starting
with an $A$-dg module $M$, we are going to compute first the dgc algebra minimal
model  of $A$:
$$
\rho \colon A' \stackrel{\simeq}{\TO} A
$$ (i.e., the classical Sullivan minimal model). Then, by means of the dgc algebra
{\em quis} $\rho$, we will make $M$ an $A'$-dg module by defining the product with
elements of $A'$ by
$$ a \cdot m = \rho^* (a) \cdot m,
$$ where the product in the right-hand side is the product of $M$ as an $A$-module.
Let us note $\rho^* (M)$ for the dg module $M$ with this structure of $A'$-module.
Finally, we will compute the minimal model of  $\rho^* (M)$ as an $A$'-dg-module: 
$$
\phii \colon M' \stackrel{\simeq}{\TO} \rho^*(M).
$$ In other words, we will compute the minimal model of the couple $(A,M)$ in the
category {\bf DGM} of {\em modules  over all algebras}. The objects of this
category are couples like $(A,M)$. Morphisms are also couples
$$ (f,\phii) \colon (A,M) \to (B,N)
$$ where $f \colon A \to B$ is a morphism of dgc algebras and $\phii \colon M \to N$
is an {\em dg-module $f$-morphism}; that is to say, a morphism of $A$-dg modules
$M \to f^* (N)$. In other words:
$$
\phii (a \cdot m) = f(a) \cdot \phii (m) \ \ \hbox{ for all } a\in A , \ m \in M.
$$ The algorithm we have described brings us the ``true'' minimal model
$$ (\rho,\phii) \colon (A',M') \to (A,M) $$ in the sense that the couple
$(A',M')$ is unique up to isomorphism (of {\bf DGM}) and the couple of {\em
quis} $(\rho,\phii) $ is also unique up to homotopies (of {\bf DGM}): this
follows from \cite[Theorem 3.4]{R3}, which tells us that the couple $(A,M)$
is minimal in {\bf DGM} if and only if $A$ is a minimal $\R$-dgc algebra
and $M$ is an $A$-dg minimal module.

\section{Stratifications and unfoldings}

We fix in this work an effective  smooth action $\Phi \colon \sbat
\times M \to M$ (non-trivial!). The orbit space of the
action is $B$ and $\pi \colon M \to B$ is the canonical
projection. The action  $\Phi$ induces on
$M$ a natural stratification by classifying  the 
points of $M$ according to their isotropy subgroups.  
This stratification  is invariant by the action of
$\sbat$, so the orbit space $B$ also inherits    a
stratified structure.  In this section we specify these
facts.

\prg {\bf Stratifications.} A {\em stratification} of
a  paracompact topological space $Z$ is a locally
finite collection ${\cal S}_Z$ of disjoint connected manifolds
called {\em strata}, such that

\bigskip

\begin{itemize}

\item[i)] $Z = {\displaystyle  \bigsqcup_{S \in {\cal S}_Z}
S}$.

\item[ii)] $S\cap \overline{S'}\not= \emptyset 
\Longleftrightarrow S \subset \overline{S'}$ (and
we write $S \leq S'$).

\item[iii)]  $({\cal S}_Z, \leq)$ is a partially ordered set
({\em poset}).

\item[iv)] There exists an open stratum $R$ which is the
maximum.

\end{itemize}

We shall say that $Z$ is a  {\em stratified space}.
Note that $R$, called {\em regular stratum}, is
necessarily dense. A {\em singular stratum} is an
element of ${\cal S}_Z$ different from $R$. We shall
write
${\cal S}_Z^{sing}$ for the family of singular strata and 
$\Sigma_Z
\subset Z$ for  its  union.
 The  {\em length} of $Z$, written $\len Z$,  is the
biggest  integer
$n$ for which there exists a chain $S_0 < S_1 < \cdots <
S_n$ of strata.  In particular, $\len Z =0$ if and only
if $Z$ is a manifold endowed with the stratification 
${\cal S}_Z = \{ \hbox{ connected components of }Z \}$. Notice that the
length is always finite.

A continuous map (resp. homeomorphism) $f \colon Y \to
Z$ between two stratified spaces is a morphism (resp.
{\em isomorphism}) if it sends  the strata of  $Y$  to the strata of
$Z$ smoothly (resp.
diffeomorphically). We shall write $Iso(Z)$ the group of isomorphisms
between $Z$ and itself. A morphism $f \colon Y \to Z$
induces a poset morphism
$f_{{\cal S}} \colon {\cal S}_Y \to {\cal S}_Z$ by
putting $f_{{\cal S}} (S) \supset f(S)$. We shall say
that $f$ is a {\em strict morphism} if the map $f_{{\cal S}}$
is strictly increasing.

\prgg  {\bf Examples.}  Through this work we shall find
the following kinds of stratification.

\begin{itemize}

\item[(a)] On a connected manifold $N$ we always may
consider  the
$0$-length stratification ${\cal S}_N =\{ N \}$. A
stratum $S\subset Z$ inherits from ${\cal S}_Z$  such a
stratification.

\item[(b)]  Any  open subset $W$  of  $Z$ inherits 
naturally from
${\cal S}_Z$  a stratified structure satisfying $\len W
\leq \len Z$. The stratification is   ${\cal S}_W = \{
\hbox{connected components of  } S
\cap W \ / \ S \in {\cal S}_Z\}$. Notice that the
inclusion $W
\hookrightarrow Z$ is a strict morphism.

\item[(c)]   Suppose  $Z$  compact.    On the product
$N \times cZ$, where   $ cZ$ is the {\em cone} $Z
\times [0,1[ \Big/  Z \times \{ 0\} $, we have the
stratification ${\cal S}_{N \times cZ}  =  \{ N
\times S \times ]0,1[ \ / \ S \in {\cal S}_Z \}
 \cup  \{  N \times   \{  \hbox{vertex  $\vartheta$ of
$cZ$} \} \}. $  Notice that $\len N \times cZ = \len Z +1$.
A point of $cZ$  will be denoted by
$[\! [x,t]\!]$ with $(x,t) \in Z \times [0,1]$. The
vertex
$\vartheta$ of $cZ$ is  $[\! [x,0]\!]$.

\end{itemize}

Unless otherwise stated, we assume that the  spaces
$W$, $cZ$ and $N
\times cZ$ are  endowed with the stratification
described  above. Later on, we shall show how $\Phi$
determines a natural stratification on
$M$ and $B$.

\prg {\bf Stratified pseudomanifolds.} When the strata
are assembled  conically  we find stratified
pseudomanifolds. We introduce this notion. An open
subset $W$ of a stratified space $Z$ is said to be {\em
modelled} on the stratified space $L$ if there exists an
isomorphism
$\phii \colon \R^{n}  \times cL \to W$ .  The pair
$(W,\phii)$ is said to be a {\em chart} of $Z$.  A
family of charts $\{ (W,\phii) \}$, where the family
$\{ W\}$ is a covering of $Z$, is called {\em atlas}. 

 We shall say that the stratified space $Z$ is a  {\em
stratified pseudomanifold} if  there exists a family 
$\{ L_S\}_{S \in {\cal S}_Z^{sing}}$ of stratified
pseudomanifolds such that for any point
$x\in\Sigma_Z$ we can find a chart  $(W,\phii)$   modelled
on $L_S$ with
$\phii ( 0,\vartheta) = x$, where $S$ is the stratum of
$Z$ containing
$x$. The  space $L_S$ is the {\em link} of the stratum
$S$. 

This definition makes sense because  it is made by
induction on the length of $Z$  ($\len L_S <  \len Z
$). A stratified space with $\len Z = 0$  is always a
stratified pseudomanifold. Each of the examples given
in \thesection.1.1  is a  stratified pseudomanifold
when $Z$ is a stratified pseudomanifold.   This
definition is slightly more general than that of
stratified pseudomanifold of \cite{GM} since we allow
the singular strata to have codimension 1.

\prg {\bf   Unfoldings}. The computation of  
the cohomology of a stratified pseudomanifold
$Z$ using differential forms is possible using the
controlled forms of Verona \cite{V}; but we need some extra
data on $Z$ so that these controlled forms  will make sense.
The original definition uses a system of neighborhoods of
singular strata  subjected to some compatibility conditions. A
more comprehensive and less technical alternative is presented
in \cite{S1} where a desingularisation of $Z$ is used. With
this blow-up, the controlled forms of $M$ and $B$ are more
easily related.In this work we follow this point of view.

Consider $Z$ a stratified pseudomanifold. A continuous
map  $ {\cal L}
\colon \wt{Z} \to Z$, where  $\wt{Z}$ is a (not
necessarily connected) manifold,  is an {\em unfolding}
if  the two following conditions hold: 
\begin{itemize}

 \item[1.] The restriction  $ {\cal L}_M \colon  {\cal
L}_M^{-1} (R ) 
\longrightarrow R $ is a local diffeomorphism.

\item[2.]   There exist a family of unfoldings  $ \{
{\cal L}_{L_S}
\colon  \wt{L_S} \to L_S\}_{S\in {\cal S}_Z^{sing}} $
and an atlas ${\cal A}$ of
$Z$ such that for  each chart $(U,\phii)\in {\cal A}$ 
there exists a commutative diagram

\begin{picture}(200,80)(-140,0)
\put(154,5){\makebox(0,0){$U$}}
\put(35,5){\makebox(0,0){$\R^{n}  \times c L_S$}}
\put(150,60){\makebox(0,0){${\cal L}_Z^{-1}(U)$}}
\put(35,60){\makebox(0,0){$\R^{n} \times \wt{L_S}
\times ]-1,1[$}}

\put(23,29.5){\makebox(0,0){$Q$}}
\put(140,29.5){\makebox(0,0){${\cal L}_Z$}}
\put(97.5,15){\makebox(0,0){$\phii$}}
\put(102,70){\makebox(0,0){$\wt{\phii}$}}

\put(65,5){\vector(1,0){65}}
 \put(82,60){\vector(1,0){40}}
\put(35,45){\vector(0,-1){31}}
 \put(150,45){\vector(0,-1){31}}

\end{picture}

where

\begin{itemize}

\item[a)]   $\wt{\phii}$ is a diffeomorphism  and  

\item[b)] $Q(x_1, \ldots , x_n,\wt{\zeta},t) = (x_1,
\ldots ,  x_n,[\![  {\cal L}_{L_S}(\wt{\zeta}),|t|]\!])
$.

\end{itemize} \end{itemize}

 This definition makes sense because it is made by
induction on the length of $Z$.   When $\len Z =0$ then
${\cal L}_Z$ is just a local diffeomorphism. The
restriction ${\cal L}_Z \colon  {\cal L}_S^{-1}(S)\to
S$ is a fibration with $\wt{L}_S$ as a fiber, for any
singular stratum
$S$.    

For each of the examples of 2.1.1 we have the
following unfoldings:

\begin{itemize}

\item[(a)] $\wt{N} = N$ and ${\cal L}_N =
\hbox{identity}$.

\item[(b)]  $\wt{W} = {\cal L}^{-1}_Z (W) $ and ${\cal
L}_W =
\hbox{restriction of }{\cal L}_Z $. 

\item[(c)]   $\wt{N \times cZ} = N \times \wt{Z} 
\times ]-1,1[ $  and
${\cal L}_{N \times cZ}(y,\wt{x} , t)  = (y,[\![  {\cal
L}_Z(\wt{x} ) , |t| ]\!])$. \end{itemize}

A morphism $f \colon Y \to Z$ between two stratified spaces, endowed with
unfoldings $ {\cal L}_Y \colon \wt{Y} \to Y$ and 
$ {\cal L}_Z \colon \wt{Z} \to Z$, is a {\em liftable morphism} if there
exists a smooth map $\wt{f} \colon \wt{Y} \to \wt{Z}$ with 
${\cal L}_Z \rond \wt{f} = {\cal L}_Z \rond f$. Each $\phii$ is a liftable
morphism. The inclusion $W \hookrightarrow Z$ is a liftable morphism.

From now on $Z$ denotes a stratified pseudomanifold
endowed with an unfolding  $ {\cal L}_Z \colon \wt{Z}
\to Z$.

\prg {\bf Stratifications induced by the action}. We
present the structure of stratified pseudomanifold of
$M$ and of the orbit space $B$. For technical reasons
we need to consider  {\sf only in this paragraph} a
smooth action $\Phi \colon G \times M \to M$ of  a
closed subgroup of the unit circle  $\sbat$ on a
manifold $M$.  The properties listed below follow
mainly from the Slice Theorem (see \cite{HS}).  

\medskip

$\bullet$ {\em Stratification}.  Consider on $M$ the
equivalence relation $\sim$ defined by  $ x \sim y $
if $G_x$ is equal to
$G_y$, where
 $G_z =  \{  g \in G\  / \ \Phi(g,z) = z \}$ denotes the
{\em isotropy subgroup} of a point $z \in M$.  Each of
the equivalence class of
$\sim$ is an invariant sub-manifold of $M$. The family
${\cal S}_M$ of the connected components of the
equivalence classes given by  this relation  defines a
stratification on $M$.  The family ${\cal S}_{B}=\{
\pi(S) \ / \ S \in {\cal S}_M \}$ defines a
stratification on the orbit space $B$. When $G$ is
connected the map $\pi_{{\cal S}}$ is bijective and therefore $\pi$ is a
strict morphism.

We shall write $G_S$ for the isotropy subgroup of a point
(and therefore any point) of a stratum $S$. Notice that
$G_S$ is a closed subgroup of
$G$. According to this   subgroup there are three types
of strata:  {\em regular stratum} ($G_S$ is $1$),  {\em
exceptional stratum} ($G_S$  is finite different from
1) and  {\em fixed stratum} ($G_S$ is $\sbat$).  Notice
that the restriction of the canonical projection $\pi
\colon M \to B$ to $S$ is a principal fibration over
$\pi(S)$ with fiber $G/G_S$. 
We shall write $F$ for the union of fixed strata. We shall identify $F \subset
M$ with $\pi (F) \subset B$ by $\pi$.

\smallskip

$\bullet$ {\em Links}.  For any singular stratum
$S\subset M$ fix a point $x$ on it and put $\S^{n_S}$
the unit sphere of a slice transversal to the stratum
$S$ at $x$. The action  $\Phi$ induces the orthogonal
action  $\Phi_S\colon G_S \times \S^{n_S} \to
\S^{n_S}$; this action has not fixed points ({\em
almost free action}).   Notice  $n_S$ is necessarily
even for a fixed stratum. The link  of $S$ is the
sphere $\S^{n_S}$ endowed with the stratification
induced by
$\Phi_S$. The link of $\pi(S)$ is the quotient space 
$\S^{n_S}/G_S$. Notice that this link is homologically 
a sphere or a real projective space when $S$ is an
exceptional stratum \cite{Br1} and a complex projective
space when $S$ is a fixed stratum \cite{NR}.

\smallskip

$\bullet$ {\em Unfoldings}. It is proven in \cite{HS}
that $M$ possesses an equivariant unfolding ${\cal L}_M
\colon \wt{M} \to M$ (relative to a free smooth
action $\wt\Phi \colon G \times
\wt{M} \to \wt{M}$) in such a way that the induced map 
${\cal L}_B
\colon \wt{M}/G \to B$ is an unfolding of $B$.
Moreover, if
$\wt\pi \colon \wt{M} \to \wt{M}/G$ is the canonical
projection, we have  ${\cal L}_B \rond \wt\pi  = \pi
\rond{\cal L}_M$. So, the morphism $\pi$ is  liftable.

 \section{Controlled forms} Controlled forms were introduced by Verona to
compute the cohomology of a stratified pseudomanifold $Z$ using differential
forms
\cite{V}. We present this notion in this paragraph, following the approach of
\cite{S1}. 

 \prg{\bf Definitions}. A differential form $\omega$ on the regular stratum $R$ of
$Z$ is said to be {\em liftable} if there exists a differential form $\wt{\omega}$ on
$\wt{Z}$, called the {\em lifting} of $\omega$, verifying:  
$ 
\wt{\omega} = {\cal L}_Z^*\omega$ on  ${\cal L}_Z^{-1}(R)$. 
 By density the lifting is unique. The differential form
$\omega$ can be tangential or transversal to the strata; in the first case we get
controlled forms and in the second case we get perverse forms.

A liftable form $\omega$ is a {\em controlled form} if it induces a differential form
$\om_S$ on each singular stratum $S$, that is, 
$\wt{\om}\big|_{{\cal L}^{-1}_Z (S)} = {\cal L}^*_Z
\om_S$. So,  we can see $\omega$ as the family of differential forms 
$\{ \om_S\in  \hiru{\Om}{}{S}\} _{S\in {\cal S}_Z}$.

We shall write $\hiruv{\Om}{}{Z}$ the {\em complex of controlled forms} (or the
{\em deRham-Verona complex}). This subcomplex of the deRham complex
$\hiru{\Om}{}{R}$ is in fact a dgc algebra. To see that, notice that if the differential forms
$\omega$ and $\eta$ are controlled, then the differential forms
$\omega + \eta$,
$\omega \wedge \eta$ and
$d\omega$  are also controlled since, for each stratum $S$, they verify: 
$\wt{\omega + \eta}\big|_{{\cal L}^{-1}_Z (S)} = 
\wt{\omega}\big|_{{\cal L}^{-1}_Z(S)} +
\wt{\eta}\big|_{{\cal L}^{-1}_Z (S)} = 
\om_S+ \eta_S$,
$\wt{\omega \wedge \eta}\big|_{{\cal L}^{-1}_Z (S)} = 
\wt{\omega}\big|_{{\cal L}^{-1}_Z (S)}\wedge
\wt{\eta}\big|_{{\cal L}^{-1}_Z (S)} =
\om_S \wedge \eta_S$,  and
$\wt{d\omega}\big|_{{\cal L}^{-1}_Z (S)} =  d\wt{\omega}\big|_{{\cal L}^{-1}_Z (S)} =
d\om_S$. 

For a stratum $S$,  we have the restriction operator  
$R_S \colon \hiruv{\Om}{}{Z} \to
\hiru{\Om}{}{S}$, defined by $R_S (\om) =
\om_S$, which is a dgc algebra operator and therefore endows 
$\hiru{\Om}{}{S}$ with a structure of $ \hiruv{\Om}{}{Z}$-dg module.

\bP 
\label{onto}
When $S$ is closed the restriction operator $R_S \colon
\hiruv{\Om}{}{Z}
\to
\hiru{\Om}{}{S}$ is onto.
\eP

\pro Fix $(U,\phi)$ a chart of ${\cal A}$. Consider $\alpha \colon \R^n \to \R$ and
$\beta \colon ]-1,1[ \to \R$ two smooth maps taking the value 1 on a neighborhood
of 0. The map
$f \colon Z-\Sigma \to \R$ defined from 
$f \phi (x_1,\ldots,x_n,[\zeta,t]) = \alpha (x_1,\ldots,x_n) \cdot \beta (t)$  is a
controlled form. In fact, its lifting is the smooth map $\wt{f} \colon \wt{Z} \to \R$
defined from 
$\wt{f} \wt\phi (x_1,\ldots,x_n,[\wt\zeta,t]) = \alpha (x_1,\ldots,x_n) \cdot \beta
(|t|)$ which is constant on the fibers of ${\cal L}_Z$. A standard argument shows that
there exists a partition of unity subordinated to ${\cal A}$.

We reduce the problem to show that 
$R_S \colon \hiruv{\Om}{}{U} \to \hiru{\Om}{}{U \cap S}$ is onto. Since $S$ is
closed then $U\cap S$ is the lowest stratum of $U$ and therefore the question
becomes: is the restriction operator
$R_{\R^n} \colon \hiruv{\Om}{}{\R^n \times cL_S} \to
\hiru{\Om}{}{\R^n}$  onto? And the answer is clearly yes. \qed

\prg {\bf Relative controlled forms.} Consider $Y$ a union of strata of $Z$. A {\em
relative controlled form} on
$(Z,Y)$ is a controlled form on $Z$ vanishing on $Y$, that is,
$\om_S\equiv 0$ for each stratum $S\subset Y$. We shall write $\hiruv{\Om}{}{Z,Y}$
the complex of relative controlled forms, which is a dgc algebra. If we consider the
restriction operator 
$R_Y = {\displaystyle \prod_{S\subset Y}} R_S
\colon   \hiruv{\Om}{}{Z} \to
\hiru{\Om}{}{Y={\displaystyle \bigsqcup_{S\subset Y} S}}= {\displaystyle
\prod_{S\subset Y}} 
\hiru{\Om}{}{S}$ then we can write
$\hiruv{\Om}{}{Z,Y} = \Ker R_Y$.

The wedge product 
$
\wedge \colon \hiruv{\Om}{}{Z} \times
\hiruv{\Om}{}{Z,Y} -\!\!\!-\!-\!\!\!\to
\hiruv{\Om}{}{Z,Y}
$ endows $\hiruv{\Om}{}{Z,Y}$ with a structure of 
$\hiruv{\Om}{}{Z}$-dg module. The natural inclusion
$
 \hiruv{\Om}{}{Z,Y} \hookrightarrow
\hiruv{\Om}{}{Z}
$ is a morphism in the category of
$\hiruv{\Om}{}{Z}$-dg modules.

\bigskip

\prg {\bf Controlled model}. The deRham-Verona complex $\hiruv{\Om}{}{Z}$ 
depends on the unfolding chosen, but for a  stratified space its cohomology does
not:
$\hiru{H}{}{\hiruv{\Om}{\cdot}{Z}}\cong\hiru{H}{}{Z}$, the singular cohomology with
real coefficients  (cf.\cite{V},\cite{S1} where {\em controlled} becomes {\em zero
perversity}).  But we have a stronger result at the level of dgc algebra minimal
models. Recall that the dgc algebra minimal model $\dos{A}{Z}$ of $Z$ is just the
dgc algebra minimal model of $A_{\R} (Z)$, the dgc algebra of polynomial
differential forms on the simplicial set
$\underline{Sing}(Z)$ of singular simplices of $Z$. The dgc algebra $\hiruv{\Om}{}{Z}$ is
easier to handle than $A_{\R}(Z)$, but we need to know that their dgc algebra
minimal models are the same.  This follows immediately from the following
Theorem (cf. \cite{Hal}).
    
\bT  
\label{abs}  Let $Z$ be a stratified space. Then there exist two {\rm quis} of dgc
algebras
$$
\hiruv{\Om}{}{Z} \stackrel{\rho_1}{\TO}  \cdot
\stackrel{\rho_2}{\longleftarrow}  A_{\R} (Z).
$$
\eT

\pro See Appendix. \qed

\prg {\bf Relative controlled model}. As in the absolute case (cf. 
Theorem \ref{abs}), the
model of a morphism $f \colon Z' \to Z$ between two stratified spaces can be
computed, under certain conditions, using controlled forms instead of polynomial
forms. These conditions involve the unfolding: 

\begin{itemize} 
\item[{[P1]}] $f$ preserves controlled forms: 
$f^* \omega \in  \hiruv{\Om}{}{Z'}$ for any $\omega \in  \hiruv{\Om}{}{Z}$
 \item[{[P2]}] $f$ preserves liftable simplices (see Appendix for the 
 exact definition).
\end{itemize} 
We shall say that $f$ satisfying the two conditions is {\em good}.

For good morphisms at least, we have a relative version of Theorem \ref{abs}.
\bT
\label{diagrama}
Let $f \colon Z' \to Z$ be a good morphism of stratified spaces. Then there exists
a commutative diagram of dgc algebra morphisms

$$ 
\begin{picture}(150,60)(0,0)
\put(-20,50){\makebox(0,0){$\hiruv{\Omega}{}{Z'}$}}
\put(50,50){\makebox(0,0){$\cdot$}}
\put(50,0){\makebox(0,0){$\cdot$}}
\put(-20,0){\makebox(0,0){$\hiruv{\Omega}{}{Z}$}}
\put(120,50){\makebox(0,0){$A_{\R}(Z')$}}
\put(120,0){\makebox(0,0){$A_{\R}(Z)$}}

\put(-22,10){\vector(0,1){30}} 
\put(50,10){\vector(0,1){30}}
\put(120,10){\vector(0,1){30}}
\put(0,0){\vector(1,0){36}} 
\put(0,50){\vector(1,0){36}}
\put(100,0){\vector(-1,0){36}} 
\put(100,50){\vector(-1,0){36}}

\put(18,8){\makebox(0,0){$\rho_1$}}
\put(18,58){\makebox(0,0){$\rho'_1$}}
\put(82,8){\makebox(0,0){$\rho_2$}}
\put(82,58){\makebox(0,0){$\rho'_2$}}

\put(-12,25){\makebox(0,0){$f^*$}} 
%\put(42,25){\makebox(0,0){$f^*$}}
\put(112,25){\makebox(0,0){$f^*$}} 

\end{picture} 
$$

in which the horizontal arrows are {\rm quis}.  
\eT

\pro See Appendix. \qed

The two examples of good morphisms used in this work are described in the
following proposition.

\bP
\label{good}

The projection $\pi \colon M \to B$ and the inclusion $\iota \colon F
\hookrightarrow B$ are good morphisms.
\eP

\pro See Appendix. \qed

\section{Models for $M$ and $F$} This section is devoted to  constructing the
minimal models of $M$ and $F$ in the category of differential graduate models
over the dgc algebra minimal model $\dos{A}{B}$ of $B$. As in other
contexts (Euler-PoincarŽ characteristic, Localization Theorem, rational homotopy
theory,...), these two models are intimately related, in fact, they are free 
$\dos{A}{B}$-graded modules over the same (up to a shift)  graded vectorial space.

\prg {\bf Model of a stratified space}.  Let $Z$ be a stratified space. As we showed
in Theorem \ref{abs}, its  minimal model as a dgc algebra can be computed from
$\hiruv{\Om}{}{Z}$. So let 
$$
\rho_Z \colon \dos{A}{Z}
\stackrel{\simeq}{\TO} \hiruv{\Om}{}{Z}
$$ be a {\em quis} of dgc algebras, with $\dos{A}{Z}$ a minimal one. Next, let 
$ f \colon Z' \to Z$ be a good morphism. Then we can endow $ \hiruv{\Om}{}{Z'}$ with
a structure of $\dos{A}{Z}$-dg module by means of the composition
$$
\dos{A}{Z} \stackrel{\rho_Z}{\TO} \hiruv{\Om}{}{Z}
\stackrel{f^*}{\TO} \hiruv{\Om}{}{Z'}.
$$ Since we will always consider this structure of module in $\hiruv{\Om}{}{Z'}$ we
shall not write 
$
\rho_Z^* f^*(\hiruv{\Om}{}{Z'})
$ but simply $\hiruv{\Om}{}{Z'}$.  In the same way, we also consider 
$A_{\R} (Z')$ as an $\dos{A}{Z}$-dg module.

\bP 
The  minimal models of $\hiruv{\Om}{}{Z'}$ and $A_{\R}
(Z')$ as $\dos{A}{Z}$-dg modules are isomorphic.
\eP
\pro It follows from Theorem \ref{diagrama}. The commutativity of the diagram
given by this result means that $\rho'_1$ and $\rho'_2$ are {\em quis} of
$\dos{A}{Z}$-dg modules.
\qed

Let us denote by $\dos{M}{Z'}$ this  minimal model as $\dos{A}{Z}$-dg module. It
obviously depends on $f$ but not on the several choices we have made in this
construction: the dgc algebra $\dos{A}{Z}$, the {\em quis} $\rho_Z$ and the
$\dos{A}{Z}$-dg module minimal model of $\hiruv{\Om}{}{Z'}$ (or $A_{\R} (Z')$).
Despite  all these choices, the
$\dos{A}{Z}$-dg module minimal model $\dos{M}{Z'}$ is unique up to isomorphism by
section  1.3.

\prg {\bf Minimal model of $F$.} The fixed point set $F$ plugs into the category of
$\dos{A}{B}$-dg modules through the natural inclusion $\iota \colon F
\hookrightarrow B$, which is a good morphism. We have already seen in
section 3.2  that $\hiruv{\Om}{}{B,F}$ is an
$\hiruv{\Om}{}{B}$-dg module and therefore an $\dos{A}{B}$-dg module. We shall
write $\dos{M}{B,F}$ for the {\em relative minimal model} of $(B,F)$, that is, the
 minimal model of
$\hiruv{\Om}{}{B,F}$ as  $\dos{A}{B}$-dg module.  Notice that the inclusion $i \colon
\hiruv{\Om}{}{B,F}
\hookrightarrow
\hiruv{\Om}{}{B}$ is an $\dos{A}{B}$-dg module morphism. We shall write
$i'
\colon
\dos{M}{B,F}
\to
\dos{A}{B}$  any of its models (see Corollary \ref{loose}). The degree of $i$ and $i'$
is 0. 

\bP
\begin{equation}
\label{F}
\dos{M}{F} = \dos{A}{B} \oplus_{i'} \dos{M}{B,F}.
\end{equation}
\eP
\pro Consider the exact sequence
$ 0 \to \hiruv{\Om}{}{B,F} \stackrel{i}{\TO} \hiruv{\Om}{}{B} 
\stackrel{R_F}{\TO} \hiru{\Om}{}{F} \to 0
$  (cf. Proposition \ref{onto}).  Then, by Remark 1.1.1 we have a $\dos{A}{B}$-dg
module {\em quis} between 
$ \hiruv{\Om}{}{B} \oplus_i \hiruv{\Om}{}{B,F}$ and $\hiru{\Om}{}{F}$. So, by
Proposition \ref{suma} we have a $\dos{A}{B}$-dg module {\em quis} between
$\dos{A}{B} \oplus_{i'} \dos{M}{B,F}$ and 
 $\hiru{\Om}{}{F}$. Since $\dos{A}{B} \oplus_{i'} \dos{M}{B,F}$  is a minimal
$\dos{A}{B}$-dg module (cf. Lemma \ref{minimal}) then  by uniqueness we get
$\dos{M}{F} =
\dos{A}{B}
\oplus_{i'} \dos{M}{B,F}$ (cf. Corollary \ref{uno}). \qed

\bC
\label{toppybis}  Suppose $B$ of finite type and simply connected, then 
$\dos{M}{B,F}$ is the free $\dos{A}{B}$-graded module generated  by $\hirut{H}{*
-1}{Y_\iota}$.
\eC
\pro From section 1.3.9 we know that $\dos{M}{F}$ is a free
$\dos{A}{B}$-graded module over $ \hiru{H}{*}{Y_\iota}$ and therefore
$\dos{M}{F} = \dos{A}{B}
\oplus_h \left[\dos{A}{B} \otimes \hirut{H}{}{Y_\iota} \right]$ for some
$h$ of degree 1. From the above Proposition we get the result wanted.
\qed

\prgg {\bf Remarks}

$\bullet$ $M$ of finite type $\Rightarrow$ $B$ of finite type. From \cite{Br2} we
know that if $M$ is of finite type then $F$ and $(B,F)$ are of finite type. Using the
long exact sequence associated to $(B,F)$  one gets that $B$ is also of finite type.

\medskip

$\bullet$ $M$ simply connected  $\Rightarrow$ $B$ simply connected.  Considering
twisted neighborhoods of orbits (cf. \cite{Br2}) one easily checks that a loop on
$B$ lifts in a path on $M$. Since the orbits of $M$ are connected, we are done.

\prg {\bf Minimal model of \mbox{\boldmath $M$}}.  The fundamental vector field 
$X$ of the action is  defined by 
$X(x) = T_1 \Phi_x (1)$, where $x \in M$ and $ \Phi_x \colon
\sbat \to M$ is given by $\Phi_x (g) = \Phi(g,x)$. Since this vector field does not
vanish on $R$  we can consider the dual form
$\chii \in \hiru{\Om}{1}{R}$, relatively a riemannian metric
$\mu$ on $R$. When this metric is {\em good} (cf. \cite{HS}) the derivative
$d\chii$ is a basic form relatively to the projection  $\pi \colon R \to \pi(R)$. So,
there exists a differential form $e\in\hiru{\Om}{2}{\pi(R)}$ verifying
$d\chii =\pi^*e$. Both differential forms, $\chii$ and $e$, are liftable. We shall say
that $e$ is an {\em Euler form}.  They are not controlled forms because their
restrictions to the links of fixed strata do not  necessarily vanish. But the maps
$
\chii \colon \hiruv{\Om}{}{M,F} \to \hiruv{\Om}{}{M}
$ and 
$ e \colon \hiruv{\Om}{}{B,F} \to \hiruv{\Om}{}{B},
$ given by $\gamma \mapsto \chii \wedge \gamma$ and $\om \mapsto e
\wedge
\om$, are well-defined $\dos{A}{B}$-dg module morphisms. We shall write 
$ e' \colon \dos{M}{B,F} \to \dos{A}{B}
$ a model of $e$.  Notice that  the degree of $e$  and $e'$ is  2. 

\medskip

The main result of this work is 

\medskip

\bT
\label{nagusia}
\begin{equation}
\label{M}
\dos{M}{M} = \dos{A}{B} \oplus_{e'} \dos{M}{B,F}.
\end{equation}
\eT
\pro Put $I\hiruv{\Om}{}{M} = \{ \om \in \hiruv{\Om}{}{M} / L_X \om =0\}$ the
complex of invariant controlled differential forms. We have seen in
\cite{HS} that the inclusion 
$I\hiruv{\Om}{}{M} \hookrightarrow \hiruv{\Om}{}{M}$ is a dgc algebra {\em quis}.
We endow
$I\hiruv{\Om}{}{M}$ with the natural structure of $\dos{A}{B}$-dg module by means
of the composition
$$
\dos{A}{B} \stackrel{\rho_B}{\TO} \hiruv{\Om}{}{B}
\stackrel{\pi^*}{\TO} \hiruv{\Om}{}{M},
$$ which is well defined since $\pi^* (\hiruv{\Om}{}{B}) \subset I\hiruv{\Om}{}{M}$.
The inclusion $I\hiruv{\Om}{}{M} \hookrightarrow
\hiruv{\Om}{}{M}$ is now a  {\em quis} of $\dos{A}{B}$-dg modules . Each invariant
differential form
$\om$ is written uniquely as $\om = \pi^* \alpha +
\chii \wedge  \pi^*\beta$; when $\om$ is controlled then $\alpha \in
\hiruv{\Om}{}{B}$ and
$\beta
\in \hiruv{\Om}{}{B,F}$ because $X$ is tangent to the links of fixed strata. So, the
operator 
\begin{equation}
\label{magique}
\Delta \colon \hiruv{\Om}{}{B} \oplus_e \hiruv{\Om}{}{B,F} \TO I\hiruv{\Om}{}{M},
\end{equation} given by $\Delta (\alpha, \beta) =  \pi^* \alpha + 
\chii \wedge \pi^*\beta$ is an $\dos{A}{B}$-dg module isomorphism and therefore 
$\dos{A}{B} \oplus_{e'} \dos{M}{B,F}$ is a
  model of $M$  as $\dos{A}{B}$-dg module (cf. Proposition \ref{suma}). Minimality
again follows from Lemma \ref{minimal}.
\qed

\prgg {\bf Remarks}.

a) Formul¾ \refp{F} and \refp{M} show that $\dos{M}{F}$ and  
$\dos{M}{M}$ are free $\dos{A}{B}$-graded modules over the same (up to a shift by
2) basis.

\smallskip

b) This theorem contains the classic result saying that, when the action is {\em
almost free} (that is, $F=\emptyset$),  the dgc algebra minimal model of
$M$ is $\dos{A}{B} \otimes \Lambda (x)$ with $\deg x=1$ (cf. \cite{GMo}).  Let us see
that.

In this case $\dos{M}{M} = \dos{A}{B} \oplus_{e'} \dos{A}{B}$ and $e'$ 
is the multiplication by a certain $e \in \dos{A}{B}$ of degree 2. 
The  {\em quis} of $\dos{A}{B}$-dg modules we have constructed $\eta \colon 
\dos{M}{M} \to \hiruv{\Om}{}{M}$ verifies  $\eta (a,b) = 
a \cdot 1 + b \cdot \eta (0,1)$. 
Consider on $\dos{M}{M}$ the product given by $(a,b) \cdot (a',b') = (a\cdot a', 
 a \cdot b' + (-1)^{\deg a'} b \cdot a')$. A 
straightforward calculation gives that $\eta$ is a {\em quis}  of dgc algebras. 
But the two dgc algebras 
 $\dos{M}{M}$
and  $\dos{A}{B} \otimes \Lambda (x)$
 are quasi isomorphic by using 
$
(a,b)\mapsto
(a \otimes 1 + b \otimes x) .
$

\medskip

We establish now some consequences of these results. 

\prg {\bf PoincarŽ polynomial}. Given a topological space
$X$ we shall write $P_X$ its PoincarŽ polynomial, that is, 
$ P_X (t) = {\displaystyle \sum_{n\geq 0}} \dim \hiru{H}{n}{X} \cdot t^n.
$
\bC
\label{polyn} Suppose $B$ of finite type and simply connected, then 
$ 1 - P_{Y_\pi} = t^2(1-P_{Y_\iota}).
$
\eC
\pro We have seen  in Proposition \ref{toppybis} that $\dos{M}{B,F}$ is a free
$\dos{A}{B}$-graded module over $\hirut{H}{*-1}{Y_\iota}$. Applying the same
method to $\pi$ we get  that $\dos{M}{B,F}$ is a free
$\dos{A}{B}$-graded module over $\hirut{H}{*+1}{Y_\pi}$. So, $\dim
\hiru{H}{0}{Y_\pi} =1$, $\dim \hiru{H}{1}{Y_\pi} =0$, $\dim
\hiru{H}{2}{Y_\pi} =\dim \hiru{H}{0}{Y_\iota} -1$ and $\dim \hiru{H}{n}{Y_\pi} =\dim
\hiru{H}{n-2}{Y_\iota}$ for $n\geq 2$. This gives the result.
\qed

\prg {\bf Vanishing of  the Euler class}. Actions with vanishing Euler class $[e]
\in \hiru{H}{2}{B-F}$ have a particular status (cf. \cite{HS} for a geometrical
interpretation). We show in the sequel how
$\dos{M}{M}$ contains information about the Sullivan minimal model of $M$ in
this case. 

When the Euler class vanishes we can choose a convenient riemannian metric on
$M$ so that $e$ itself vanishes (cf. \cite{HS}). Thus $e'=0$. The minimal model 
$\dos{M}{M}$ is of the form $\dos{A}{B} \otimes E$ with $E^{^0} =\R$ and  $d E
\subset \dos{A}{B}
\otimes E^{^+}$. It supports a dgc algebra structure by putting on $E$ the trivial
product:
$1\cdot v = v$ if $v \in E$ and $E^{^+} \cdot E^{^+} =0$. This dgc algebra structure
shall be said {\em na•ve}. It contains the following information about
$M$.

\bC   
If the Euler class vanishes then the na•ve dgc algebra structure of
$\dos{M}{M}$  has the same real homotopy type of $M$. Moreover,
$\pi^{*}_\psi (B)$ injects into $\pi^{*}_\psi (M)$.
\eC
\pro  The operator $\Delta \colon \hiruv{\Om}{}{B} \oplus_0
\hiruv{\Om}{}{B,F}
\TO \hiruv{\Om}{}{M}$ is a {\em quis}  of  $\dos{A}{B}$-dg modules which becomes a
 {\em quis}  of dgc algebras when considering on the source the following product:
\begin{equation}
\label{prod} (\alpha,\beta)\cdot (\alpha',\beta') = (\alpha \cdot \alpha' , (-1)^{\deg
\alpha}
\alpha\cdot \beta' +(-1)^{\deg \alpha' \cdot \deg \beta} \alpha'
\cdot \beta).
\end{equation}   
This dgc algebra contains the real homotopy type of $M$.

Put $\rho_{(B,F)} \colon \dos{M}{B,F} \to \hiruv{\Om}{}{B,F}$ the relative minimal
model of $(B,F)$. The operator
$$ \rho_B \oplus \rho_{(B,F)} \colon \dos{M}{M} = \dos{A}{B} \oplus_0
\dos{M}{B,F} \TO  \hiruv{\Om}{}{B} \oplus_0 \hiruv{\Om}{}{B,F}
$$  
is a
{\em quis} of $\dos{A}{B}$-dg modules which becomes a  {\em quis} of dgc algebras
when considering  the product \refp{prod} in both terms. This dgc algebra contains
the real homotopy type of $M$. Notice that the dgc algebra structure on
$\dos{M}{M}$ given by \refp{prod} is just the na•ve structure. This gives the first
part of the Corollary.

\medskip

 We shall  write $\dos{A}{B} = \Lambda Y$, with differential $\partial$, and
$\dos{M}{M} = \Lambda Y \otimes E$, with differential $d$. We have $d_{|
\Lambda Y} = \partial$ and $d E \subset \Lambda Y \otimes E^{^+}$. This last
property allows us to construct a
 KS-extension

\bigskip

$$ 
\begin{picture}(100,40)(0,0) 

\put(0,40){\makebox(0,0){$(\Lambda Y,\partial)$}} 
\put(100,40){\makebox(0,0){$(\Lambda Y \otimes E,d)$}}
\put(100,0){\makebox(0,0){$(\Lambda Y \otimes \Lambda X, \delta)$}}

\put(20,40){\vector(1,0){45}}
\put(15,30){\vector(2,-1){45}} 
\put(100,10){\vector(0,1){20}}

\put(110,20){\makebox(0,0){$\phii$}}

\end{picture}                                           
 $$

\bigskip

\nt verifying $\phii_{|\Lambda Y} = \id_{\Lambda Y}$ and  $\phii (X) \subset
\Lambda Y  \otimes E^{^+}$.  Recall that $\pi_\psi^*(B) =
\hiru{H}{}{Y,\partial_0}$ and  $\pi_\psi^*(M) =
\hiru{H}{}{Y \oplus X,\delta_0}$, where $\partial_0$ and
$\delta_0$ are the linear part of $\partial$ and $\delta$ respectively. 

Since
$(\Lambda Y,\partial)$ is minimal we just have $\partial_0 = 0$ and therefore 
$\pi_\psi^*(B) = Y$. On the other hand, the composition
$
\delta_0 = X \stackrel{\delta}{\to} \Lambda  Y \otimes \Lambda X
\stackrel{\proj}{\TO} \Lambda Y \stackrel{\proj}{\TO} Y 
$ vanishes: if $x\in X$ then $ \phii (\delta x) = d (\phii x) \in d (\Lambda Y \otimes
E^{^+}) \subset \Lambda Y \otimes E^{^+}$. So, $\pi_\psi^* (M) = Y
\oplus
\hiru{H}{}{X,\delta_0}$ and the proof is finished  since $\pi^*$ becomes the
inclusion 
$Y \hookrightarrow Y \oplus \hiru{H}{}{X,\delta_0}$.
\qed

\prgg {\bf Remarks}. 

a)  When $B$ is also contractible then $\dos{M}{M} = E$ is just a dgc algebra with
trivial product and, by minimality, with zero derivative. In other words,
$M$ is a wedge of spheres.

\smallskip

b) The Gysin sequence associated to \refp{M}  implies that the cohomology of $B$
injects into the cohomology of $M$; in fact, we have the short exact sequence
$$ 0 \to \hiru{H}{*}{B}  \to \hiru{H}{*}{M} \to \hiru {H}{*-1}{B,F} \to 0.
$$  The last statement  of the Corollary implies that, when $M$ is of finite type and
simply connected, we have the following short exact sequence 
$$  0 \to \pi^* (B) \otimes \R  \to \pi^* (M) \otimes \R \to \pi^*(Y_{\pi}) \otimes
\R \to 0.
$$

\smallskip

c) The na•ve structure of $\dos{M}{M}$ appears when $e'=0$, but the previous result
needs the vanishing of the Euler class itself as is shown in the following example. 
Consider the action $\Phi \colon \sbat  \times \C\P^n \to \C\P^n$ given by 
$ z \cdot [z_0, z_1, \ldots , z_n] = [z_0, z\cdot z_1, \ldots , z\cdot z_n]
$ (in homogeneous coordinates). Here the fixed point set is $F = \C\P^0  \cup
\C\P^{n-1}
$ and the orbit space $B$ is the closed cone over
$\C\P^{n-1}$. So, $B$ is acyclic and $e' =0$. The Euler class does not vanish since
it generates
$
\hiru{H}{2}{B-F} = \hiru{H}{2}{\C\P^{n-1} \times ]0,1[} =  \R$.  The minimal model we
have computed in Theorem \ref{nagusia} is
$
\dos{M}{\C\P^n} = \R \oplus_0 \hirut{H}{*-1}{\C\P^0  \cup \C\P^{n-1}}.
$ Considering on it the na•ve structure we get
$
\dos{M}{\C\P^n} =  \hiru{H}{*}{\S^2 \vee \S^4 \vee \cdots \vee \S^{2n}}
$
 as dgc algebras. But clearly this dgc algebra does not contain the real
homotopy type of $\C\P^n$.

\prg {\bf Cohomological dimension}. Write $\dimc (X)$ for the cohomological
dimension of the topological space $X$, that is, $\dimc (X) = \sup \{ n
\in \N  \  / \ \hiru{H}{n}{X} \not = 0\}$.

\bC    Under the hypothesis of Corollary \ref{polyn}, if $\dimc( B)$ and $\dimc
(Y_\pi) $ are finite  then
$\dimc (M) = \dimc (F)$ or $\dimc (F) +2$.
\eC

\pro  From Corollary \ref{polyn} we get 
$\dimc (Y_\pi) = \dimc (Y_\iota)=0$ or
$\dimc (Y_\pi) = \dimc (Y_\iota)+2$. Now considering  the homotopy fibrations
associated to $\pi$ and
$\iota$ we get $\dimc (M) = \dimc (B)+\dimc (Y_\pi)$ and 
 $\dimc (F) = \dimc (B)+\dimc (Y_\iota)$ and then we get the result.
\qed

\medskip

We find examples of this situation when $M= \C^n $, where $\sbat$ 
acts by complex multiplication, and 
$M= \S^{n+2} = \S^1 * \S^n$, where $\sbat$ acts by multiplication on the first
factor and trivially on the second factor. When $M$ is compact and oriented the
condition  $\dimc (M) =
\dimc (F) =0$ does not occur and the condition  $\dimc (M) =
\dimc (F) +2$ is equivalent to saying that $F$ possesses a connected component of
codimension 2. This is also equivalent to the fact that $B$ has a boundary.  So,
under the conditions of  Corollary \ref{polyn}, if
$\dimc (B )< \infty$  and $B$
 has no  boundary then $\dimc (Y_\pi) = \dimc (Y_\iota) = \infty$.

\prg{\bf Equivariant cohomology}. 

 The {\em equivariant cohomology} of $M$ is the cohomology of the quotient space
$M_{\sbat} = M \times_{\sbat}
\S^\infty$, written $\hirue{H}{}{M}$. The natural projection $p \colon M_{\sbat}
\to B$ induces a natural structure of
$\dos{A}{B}$-dg module on $M_{\sbat}$. Here we compute  the {\em equivariant
minimal model} of $M$, that is,
$\dos{M_{\sbat}}{M} = \dos{M}{M_{\sbat}}$.

We shall write $\Lambda(e)$ the polynomial algebra generated by an element $e$ of
degree 2. The trivial
$\dos{A}{B}$-dg module structure will be considered on it. The main result in this
framework is

\bT If the fixed point set $F$ is not empty, then
\begin{equation}
\label{E}
\dos{M_{\sbat}}{M} = \left[ \dos{A}{B} \otimes \Lambda(e) \right]
 \oplus_{q'}
\left[\dos{M}{B,F}\otimes \Lambda(e) \right],
\end{equation} where $q'(b\otimes e^n) =e'(b) \otimes e^n + i'(b) \otimes e^{n+1}$.
\eT

\pro The equivariant cohomology can be computed using the complex 
$
\hirue{\Om}{}{M} = I\hiru{\Om}{}{M}\otimes \Lambda(e),
$  endowed with the derivative  $d(\om \otimes e^n) = (d\om) \otimes e^n + (i_X
\om) \otimes e^{n+1}
$. Here $i_X$ denotes the contraction by $X$. Proceeding as in
\cite[p.213]{S1} one shows that the two dgc algebras $I\hiru{\Om}{}{M}$ and
$I\hiruv{\Om}{}{M}$ are    quasi-isomorphic. Therefore the equivariant cohomology
of
$M$ is computed by using
$I\hiruv{\Om}{}{M}\otimes \Lambda (e)$ and $p$ induces the operator
$ P \colon \hiruv{\Om}{}{B} \to I\hiruv{\Om}{}{M}\otimes \Lambda (e)$ defined by
$P(\alpha) = \pi^*\alpha
\otimes 1$. Under these transformations the $\dos{A}{B}$-dg module structure on
$I\hiruv{\Om}{}{M}\otimes
\Lambda (e)$ is given by:
$$  a \cdot (\om \otimes e^n) = P( \rho_B (a)  )  \cdot (\om \otimes e^n), \ \
\hbox{with} \ \ a 
 \in \dos{A}{B}, \om \in  I\hiruv{\Om}{}{M}.
$$  
Now we compute  the minimal model of $I\hiruv{\Om}{}{M}\otimes
\Lambda (e)$ relative to  this structure.

The $\dos{A}{B}$-dg module isomorphism $\Delta \colon \hiruv{\Om}{}{B} \oplus_e
\hiruv{\Om}{}{B,F}
\TO \hiruv{I\Om}{}{M}$ induces  the $\dos{A}{B}$-dg module isomorphism 
$$
\nabla \colon  \left[ \hiruv{\Om}{}{B}  \otimes \Lambda (e)\right]
\oplus_q
\left[\hiruv{\Om}{}{B,F} \otimes \Lambda (e)\right] \TO \hiruv{I\Om}{}{M} \otimes
\Lambda (e),
$$ where $ q( \beta \otimes e^n  ) = (\beta \wedge e) \otimes  e^n+ i (\beta )
\otimes e^{n+1}.
$ Proceeding as in Theorem \ref{nagusia} we get that a model of 
$\hiruv{I\Om}{}{M} \otimes \Lambda (e)$ is just 
$\left[\dos{A}{B} \otimes \Lambda(e) \right] \oplus_{q'} \left[\dos{M}{B,F} \otimes
\Lambda(e) \right],
$ where $q'(b\otimes e^n) =e'(b) \otimes e^n + i'(b) \otimes e^{n+1}$.  This model
is minimal because of Lemma \ref{minimal}. \qed

Notice that for the almost free case ($F = \emptyset$) we have obtained the
following non-minimal $\dos{A}{B}$-dg module model: 
$ \left[ \dos{A}{B} \otimes \Lambda(e) \right]
 \oplus_{i' \otimes e}
\left[ \dos{A}{B} \otimes \Lambda(e) \right]$; the minimal one is just  $\dos{A}{B}$.

\prgg {\bf Remarks}

\smallskip

a) {\it PoincarŽ polynomial.}  When  $F \not= \emptyset$, $B$ is of finite type and
simply connected the relation between the homotopy fibers of $p$ and $\pi$ is
given by 
$ P_{Y_\pi} = (1-t^2) P_{Y_p}
$ (same proof as that of Corollary \ref{polyn}).

%\smallskip

 b) {\it Extension of scalars}. The complexes $\hirue{\Om}{}{M}$ and
$\dos{M_{\sbat}}{M}$  naturally support a structure of $\Lambda(e)$-dg module. The
$\dos{A}{B}$-dg module {\em quis} we have constructed is in fact a
$\Lambda(e)$-dg module {\em quis}. For this  structure, the extension of scalars of
$\dos{M_{\sbat}}{M}$ is just $\dos{M}{M}$, that is, 
$
\R \otimes_{\Lambda (e)} \dos{M_{\sbat}}{M} =\dos{M}{M}.
$ In other words, the model of the fiber of $M \to M_\sbat \to B$ is "the fiber of the
model" (cf. \cite{Hal}).

%\smallskip

c) {\it Vanishing Euler class}. Since $e'=0$ one gets that $\dos{M_\sbat}{M}$ is
isomorphic to 
$
\dos{A}{B} \oplus \left[ \dos{M}{F} \otimes \Lambda^+ (e) \right],
$ relative to both module structures. We conclude that the cohomology of $B$
injects into the equivariant cohomology of $M$, that is, we have in fact the
following exact sequence
$ 0 \to \hiru{H}{}{B} \to \hirue{H}{}{M} \to \hiru{H}{}{F} \otimes \Lambda^+ (e) \to 0.
$

%\smallskip

d) {\it Equivariant cohomology}.   Formula \refp{E} says that we can compute the
equivariant minimal model $\dos{M_{\sbat}}{M}$ in terms of basic data: 
$ i',e' \colon \dos{M}{B,F} \flechas \dos{A}{B}.
$ The equivariant cohomology $\hirue{H}{*}{M}$ can  also be computed  in terms of
basic data. In fact, the short exact sequence
$ 0 \to \dos{A}{B}  \otimes \Lambda (e)\to \dos{M_{\sbat}}{M} \to \dos{M}{B,F}
\otimes \Lambda (e)
\to 0
$  associated to \refp{E}  (cf. Lemma \ref{112}) gives the long exact sequence
$$
\cdots \to \left[ \hiru{H}{}{B} \otimes \Lambda (e) \right]^{^i} \to
\hirue{H}{i}{M} \to \left[ \hiru{H}{}{B,F} \otimes \Lambda (e) \right]^{^{i-1}}
 \stackrel{q^*}{\TO} \left[ \hiru{H}{}{B} \otimes \Lambda (e) \right]^{^{i+1}}
\to \cdots ,
$$ which determines $\hirue{H}{}{M}$ in terms of 
$ i^*,e^* \colon \hiru{H}{}{B,F} \flechas \hiru{H}{}{B}.
$

e) {\it Equivariantly formal spaces.} Put ${\mathfrak r} \colon M \rightarrow M_\sbat$
the inclusion given by ${\mathfrak r}(x) = \hbox{ class of } (x,1)$. The manifold $M$
is {\em equivariantly formal} if the restriction map ${\mathfrak r}^* \colon
\hirue{H}{}{M}
\to
\hiru{H}{}{M}$ is surjective (cf. \cite{Br2}, \cite{GKMP}). We can translate this
condition in terms of basic data by considering the following commutative diagram

$$ 
\begin{picture}(300,60)(0,0)
\put(-65,50){\makebox(0,0){$\cdots$}}
\put(-20,50){\makebox(0,0){$\hirue{H}{i}{M}$}}
\put(100,50){\makebox(0,0){$\left[ \hiru{H}{}{B,F} \otimes \Lambda (e)
\right]^{^{i-1}}$}}
\put(230,50){\makebox(0,0){$\left[ \hiru{H}{}{B} \otimes \Lambda
(e)\right]^{^{i+1}}$}}
\put(300,50){\makebox(0,0){$\cdots$}}

\put(-65,0){\makebox(0,0){$\cdots$}}
\put(-20,0){\makebox(0,0){$\hiru{H}{i}{M}$}}
\put(100,0){\makebox(0,0){$\hiru{H}{i-1}{B,F}$}}
\put(230,0){\makebox(0,0){$\hiru{H}{i}{B}$}}
\put(300,0){\makebox(0,0){$\cdots$}}

\put(-20,40){\vector(0,-1){30}} 
\put(100,40){\vector(0,-1){30}}
\put(230,40){\vector(0,-1){30}}

\put(-55,50){\vector(1,0){10}}
\put(10,50){\vector(1,0){30}}
\put(150,50){\vector(1,0){30}}
\put(275,50){\vector(1,0){10}}

\put(-55,0){\vector(1,0){10}}
\put(10,0){\vector(1,0){50}}
\put(145,0){\vector(1,0){60}}
\put(255,0){\vector(1,0){30}}

\put(175,10){\makebox(0,0){$e^*$}}
\put(165,60){\makebox(0,0){$q^*$}}

\put(-12,25){\makebox(0,0){${\mathfrak r}^*$}} 
\put(112,25){\makebox(0,0){${\mathfrak R}$}}
\put(242,25){\makebox(0,0){${\mathfrak R}$}} 

\end{picture} 
$$ 

\smallskip

\nt where ${\mathfrak R}(\sum [\alpha_n ] \otimes e^n) = [\alpha_n ]$. Since the
restriction 
${\mathfrak R} \colon \Coker q^* \to \Coker e^*$ is surjective then the manifold
$M$ is equivariantly formal if and only if the restriction ${\mathfrak R} \colon
\Ker q^*
\to
\Ker e^*$ is surjective. In other words, any string $[ \alpha_0 ]
\stackrel{e^*}{\rightarrow} 0$ fits into a string 
$ 0 
\stackrel{i^*}{\leftarrow} [ \alpha_n ]
\stackrel{e^*}{\rightarrow} [\beta_n ]
\stackrel{i^*}{\leftarrow} \cdots
\stackrel{e^*}{\rightarrow} [\beta_0]
\stackrel{i^*}{\leftarrow} [\alpha_0]
\stackrel{e^*}{\rightarrow}0.
$

\smallskip

f) {\it Localization Theorem.} This Theorem asserts that the restriction map
$R_F \colon M \to F$ induces an isomorphism between their localizations 
$S^{-1}\hirue{H}{}{M}$ and $S^{-1}\hirue{H}{}{F}$. We can translate this Theorem in
terms of basic data saying that the map
$$
\nabla \colon \hiru{H}{}{B,F;S^{-1}\Lambda (e) } \longrightarrow 
\hiru{H}{}{B,F;S^{-1} \Lambda (e)},
$$ defined from $\nabla ([\om]) =e \cdot [\om] + [\om \wedge e]$, is an
isomorphism. This comes from the fact that the Localization Theorem is
equivalent to the vanishing of 
$S^{-1}\hirue{H}{}{M,F}$ and from the exact sequence
$$
\cdots \to \left[ \hiru{H}{}{B,F} \otimes \Lambda (e) \right]^{^i} \to
\hirue{H}{i}{M,F} \to \left[ \hiru{H}{}{B,F} \otimes \Lambda (e) \right]^{^{i-1}}
 \stackrel{q^*}{\TO} \left[ \hiru{H}{}{B,F} \otimes \Lambda (e) \right]^{^{i+1}}
\to \cdots ,
$$ obtained proceeding as in d). In fact, for $S = \left\{ 1,e^p,e^{2p},\ldots \right\}$,
the inverse of $\nabla$ is given by \newline $\nabla^{-1} ([\om]) = {\displaystyle
\sum_{n\geq 0}} {\displaystyle \sum_{j=0}^{p-1}}  (-1)^{(n+1)p-1-j} (e^{-p})^{n+1}
e^j \cdot  [ \om \wedge  e^{p(n+1)-j-1}]$, which makes sense since the differential
form $e^m$ vanishes for large enough $m$.

\prg {\bf Semifree $\S^3$-actions}. The results developed until here for circle
actions extend directly to semifree 
$\S^3$-actions. This comes essentially from the fact that formula \refp{magique}
applies (up to a shift) here (cf. \cite{S3}). We don't give all the results but just the
main one.

Consider $\Phi \colon \S^3 \times M \to M$ a semifree smooth action of
$\S^3$ on a smooth manifold $M$. Write $F$ the submanifold of fixed points.
Semifreeness means that $\S^3$ acts freely on $M-F$. The orbit space $B$ is a
stratified pseudomanifold whose singular strata are the connected components of
$F$. The inclusion $\iota \colon F \hookrightarrow B$ induces the inclusion operator 
$ i \colon \hiruv{\Om}{}{B,F} \TO \hiruv{\Om}{}{B}
$ which is an $\dos{A}{B}$-dg module morphism of degree 0. Here the Euler form
$e$ lies on $\hiru{\Om}{4}{B-F}$ and  induces the $\dos{A}{B}$-dg module
morphism 
$ e \colon \hiruv{\Om}{}{B,F} \TO \hiruv{\Om}{}{B}
$  of degree 4. The $\dos{A}{B}$-dg module minimal models of these operators are 
$i'$ and $e'$ respectively. The main result in this framework is

\bT
$$
\dos{M}{M} = \dos{A}{B} \oplus_{e'} \dos{M}{B,F}
$$
$$
\dos{M}{F} = \dos{A}{B} \oplus_{i'} \dos{M}{B,F}
$$
\eT

\pro  Follow the path taken in the Theorem \ref{nagusia}. \qed

\prg {\bf Isometric flows}.  An {\em isometric flow} is a real smooth action
$\Phi \colon \R \times M \to M$ preserving a riemannian metric $\mu$ on the
smooth manifold $M$. The fundamental vector field 
$X$ of the action is a Killing vector field. We shall write ${\cal F}$ the singular
foliation determined by the orbits of the action. The fixed point set is a manifold
written $F$. Notice that in this  case the orbit space $B = M   / \R$ can be very wild
(even totally disconnected!). For this reason it is customary to work with ``basic
objects" (objects living on $M$ transverse to the flow and invariant by the flow)
instead of working directly  with the objects living on $B$. For example, a {\em
basic form} is a differential form on
$M$ which is transverse to the flow ($i_X\om =0$) and  invariant by the flow ($L_X
\om = 0$), a {\em basic controlled form}  is a controlled form on $M$ verifying
$i_X\om= i_X d\om =0$,
$\ldots$ We shall write
$$
\begin{array}{cl}
\hiru{\Om}{}{M/{\cal F}}& \hbox{the complex of basic forms}\\
\hiruv{\Om}{}{M/{\cal F}}& \hbox{the complex of basic  controlled forms}\\
\hiruv{\Om}{}{(M,F)/{\cal F}}& \hbox{the complex of basic relative controlled
forms.}
\end{array}
$$ When the action is periodic we have in fact a circle action and these complexes
become, up to isomorphism,
$
 \hiru{\Om}{}{B},
$
$
\hiruv{\Om}{}{B}
$ and
$
\hiruv{\Om}{}{B,F}
$ respectively. 

The three complexes above are in fact dgc algebras. The dgc algebra minimal
model of
$\hiruv{\Om}{}{M/{\cal F}}$ and $\hiru{\Om}{}{M/{\cal F}}$ are the same and they will
be denoted by
$\dos{A}{M/{\cal F}}$ (cf. \cite{S2}).  We work in the category of
$\dos{A}{M/{\cal F}}$-dg modules. The $\dos{A}{M/{\cal F}}$-dg module minimal
model of
$\hiruv{\Om}{}{(M,F)/{\cal F}}$ will be denoted by $\dos{M}{(M,F)/{\cal F}}$. 

The inclusion $\iota \colon F \hookrightarrow M$ induces the inclusion operator 
$ i \colon \hiruv{\Om}{}{(M,F)/{\cal F}} \TO \hiruv{\Om}{}{M/{\cal F}}
$ which is an $\dos{A}{M/{\cal F}}$-dg module morphism of degree 0. Here the Euler
form
$e$ lies on $\hiru{\Om}{2}{(M-F)/{\cal F}}$ and it induces the
$\dos{A}{M/{\cal F}}$-dg module morphism 
$ e \colon \hiruv{\Om}{}{(M,F)/{\cal F}} \TO \hiruv{\Om}{}{M/{\cal F}}
$  of degree 2. The $\dos{A}{M/{\cal F}}$-dg module minimal models of these
operators are 
$i'$ and $e'$ respectively.

We shall write $X= \partial /\partial t$ if there exists a diffeomorphism $M\cong  B
\times \R$ sending $X$ on a multiple of $\partial / \partial t$. The main result in
this framework is

\bT If $X \not= \partial /\partial t$ then 
$$
\dos{M}{M} = \dos{A}{M/{\cal F}} \oplus_{e'} \dos{M}{(M,F)/{\cal F}}
$$
$$
\dos{M}{F} = \dos{A}{M/{\cal F}} \oplus_{i'} \dos{M}{(M,F)/{\cal F}}
$$
\eT

\pro  If we prove that the inclusion $\hiruv{I\Om}{}{M}  \hookrightarrow
\hiruv{\Om}{}{M}$ and the restriction $\hiru{\Om}{}{M}  \to
\hiruv{\Om}{}{M}$ are dgc algebra {\em quis}, then it suffices to follow the path
taken in the proof of the Theorem \ref{nagusia}. 

An isometric flow defines a singular riemannian foliation ${\cal F}$ where the
 orbits are all closed or the closures of the orbits are all tori (cf. \cite{Mol}). In the first case the
natural projection $\pi \colon M \to M/{\cal F}$ becomes a locally trivial fibration
and, by orientability, a trivial one. So, $X=\partial  / \partial t $ and we are in the
second case. Using the Mayer-Vietoris argument we can replace $M$ by a torus
$\T$ endowed with an $\R$-linear action. The problem is then to prove that the
inclusion  $\hiruv{I\Om}{}{\T}  \hookrightarrow
\hiruv{\Om}{}{\T}$ and the restriction $\hiru{\Om}{}{\T}  \to
\hiruv{\Om}{}{\T}$ are dgc algebra {\em quis}. Since the flow is regular then 
$\hiruv{\Om}{}{\T}  = \hiru{\Om}{}{\T} $. By density, the complex $\hiruv{I\Om}{}{\T} 
=\hiru{I\Om}{}{\T} $ becomes $\hiru{\Om_{\T}}{}{\T} = \{ \om \in \hiru{\Om}{}{\T}
\hbox{ invariant by } \T \}$.  Since $\T$ is compact, we already know that the
inclusion $\hiru{\Om_{\T}}{}{\T}  \hookrightarrow
\hiru{\Om}{}{\T}$ is a  dgc algebra {\em quis} (cf. \cite{GHV}).  \qed

\medskip

In the case $X = \partial/ \partial t$ this result would imply $\hiru{H}{*}{M} =
\hiru{H}{*}{B \times\sbat}$, which is false.  We finish the work extending to
isometric flows a well known result  related to periodic actions.

\bC 
 If $X \not= \partial /\partial t$ then, for any $r\geq 0$:
$$
\hiru{H}{r-1}{(M,F)/{\cal F}} +\sum_{i=0}^{\infty} \dim \hiru{H}{r+2i}{F} \leq 
\sum_{i=0}^{\infty} \dim \hiru{H}{r+2i}{M}.
$$
\eC

\pro  From the above formul¾ we get the following  long exact sequences:
$$
\cdots \to \hiru{H}{i}{M/{\cal F}} \to \hiru{H}{i}{M} \to \hiru{H}{i-1}{(M,F)/{\cal F}}
\to \hiru{H}{i+1}{M/{\cal F}} \to \cdots,
$$
$$
\cdots \to \hiru{H}{i}{(M,F)/{\cal F}} \to \hiru{H}{i}{M/{\cal F}} \to \hiru{H}{i}{F} \to
\hiru{H}{i+1}{(M,F)/{\cal F}} \to \cdots.
$$ Since the action is free out of $F$ the above Theorem admits the  relative
version:
$
\dos{M}{M,F} = \dos{M}{(M,F)/{\cal F}} \oplus_{e'} \dos{M}{(M,F)/{\cal F}}.
$ This gives the long exact sequence 
$$
\cdots \to \hiru{H}{i}{(M,F)/{\cal F}} \to \hiru{H}{i}{M,F} \to \hiru{H}{i-1}{(M,F)/{\cal F}}
\to \hiru{H}{i+1}{(M,F)/{\cal F}} \to \cdots.
$$  Finally, by considering the long exact sequence associated to $(M,F)$ one gets
all the ingredients to proceed as in \cite[pag.161]{Br2} in order to obtain the
following Smith-Gysin sequence:
$$
\cdots \to \hiru{H}{i}{(M,F)/{\cal F}} \to \hiru{H}{i}{M} \to \hiru{H}{i-1}{(M,F)/{\cal F}}
\oplus \hiru{H}{i}{F} \to \hiru{H}{i+1}{(M,F)/{\cal F}} \to \cdots
$$ Now it suffices to follow the procedure of \cite[pag. 127]{Br2}.
\qed

When $M$ is compact the group of isometries of $(M,\mu)$ is a compact Lie group.
So, the action of $\R$ extends to an action of a torus $\T$.  For this action we have
the inequality
$$
\hiru{H}{r-1}{(M,F)/\T} +\sum_{i=0}^{\infty} \dim \hiru{H}{r+2i}{F} \leq 
\sum_{i=0}^{\infty} \dim \hiru{H}{r+2i}{M},
$$ for each $r\geq 0$ (cf. \cite{Br2}), but notice that in general
$\hiru{H}{*}{(M,F)/{\cal F}}$ and
$\hiru{H}{*}{(M,F)/\T}$ are not equal.

\section{Appendix}
In this Appendix we give the proof of 
Theorem \ref{THI} (existence of a minimal KS-extension),
Theorem \ref{THII} (uniqueness of a minimal KS-extension), 
Theorem \ref{abs} (minimal model versus controlled forms),
Theorem \ref{diagrama} (relative minimal model versus controlled 
forms) and
Proposition \ref{good} (examples of good morphisms).

\prg {\bf Proof of Theorem \ref{THI}.}
 
To begin with, put $N(0,0) =M$, $\iota _{(0,0)} = \Id_M$ and
$\rho _{(0,0)} = \phii $. Let us assume that we have already
constructed $$
\iota_{(n,0)}\colon M\longrightarrow N(n,0)
\hspace{1cm} \hbox{and}\hspace{1cm}  \rho_{(n,0)}\colon N(n,0)
\longrightarrow X $$ such that

\bigskip

$ \begin{array}{rl} (i)_{(n,0)} & \iota_{(n,0)} \hbox{ is a minimal
KS-extension},
\\[.3cm]  (ii)_{(n,0)} & \rho _{(n,0)}\iota_{(n,p)} =\phii,\hbox{ and}
\\[.3cm] (iii)_{(n,0)} & ( \rho_{(n,0)} )^{^i}_{_*}
\colon \hiru{H}{i}{N(n,0)} \to \hiru{H}{i}{X}
\hbox{ is an isomorphism for $0 \leq i \leq n-1$ and } (
\rho_{_{(n,0)}} )^{^n}_{_*}  \\[.2cm] &   \hbox{is a monomorphism}. 
\end{array} $

\bigskip

\nt Now, for $q>0$, take $ V(n,q) = \hiru{H}{n+1}{N(n,q-1),X}
$ and consider it as a homogeneous  vector space of degree $n$.
Take also a  linear section $s$ of the natural projection
$\hiru{Z}{n+1}{N(n,q-1),X}
\longrightarrow V(n,q)$. So for every $v \in V(n,q)$, we have $ sv
= (t_v, x_v)
\in
\bi{N(n,q-1)}{n+1} \oplus X^{^n} $ such that 
$$
\left( \begin{array}{c} 0 \\ 0 \end{array} \right) = \left(
\begin{array}{cc} d              &  0  \\ \rho_{(n,q-1)}  & -d
\end{array} \right) \left( \begin{array}{c} t_v  \\  x_v
\end{array}
\right). 
$$ 
Define 
$$ 
N(n,q) = N(n,q-1) \oplus  (A \otimes V(n,q)) $$ with differential $ dv = t_v
$ Define also $
\iota_{(n,q)}\colon M\longrightarrow N(n,q) $ as the composition of
$\iota_{(n,q-1)}$ with the natural inclusion of
$N(n,q-1) \hookrightarrow N(n,q)$. Finally, take $
\rho_{(n,q)}\colon N(n,q)
\longrightarrow X $ to be the morphism of $A$-dg modules
induced by 
$\rho_{(n,q-1)}$ and the linear map $f\colon V(n,q)\longrightarrow
X$ defined by
$fv = x_v$.

\medskip

\noindent Let us verify that:

\bigskip

$ \begin{array}{rl}
(i)_{(n,q)} &   \iota _{(n,q)} \hbox{ is a minimal
KS-extension because of $(i)_{(n,q-1)}$ and the definition of
$\iota _{(n,q)}$},\\[,3cm]

(ii)_{(n,q)} &   \rho _{(n,q)} \iota _{(n,q)} = \rho _{(n,q-1)} \iota
_{(n,q-1)} =\phii 
\hbox{  because of the definitions and
$(ii)_{(n,q-1)}$}, \\[,3cm]

(iii)_{(n,q)} &( \rho_{(n,q)} )^{^i}_{_*} \colon
\hiru{H}{i}{N(n,q)}
\longrightarrow \hiru{H}{i}{X} \hbox{  is an isomorphism for
$i = 0, \dots , n$}. 
\end{array}
$

\begin{itemize}
\item[]  In fact, for $i < n$, $\rho _{(n,q)}$ coincides with $\rho
_{(n,q-1)}$ and so the morphism induced in cohomology is an
isomorphism by
$(iii)_{(n,q-1)}$. In degree $n$, the natural inclusion
$N(n,q-1)
\hookrightarrow N(n,q)$ induces a monomorphism in cohomology
because new generators can only kill cocycles in degree $>n$. So
$( \rho _{(n,q)} )^{^n}_{_*}$ is a monomorphism and all we need to prove is that it is
also an epimorphism:  let
$[x] \in \hiru{H}{n}{X}$ and $[x] \not\in  \Ima (\rho_{(n,q-1)} )^{^n}_{_*}$ (otherwise
we are done). Then $x$ defines a non-zero relative cohomology class $v =  [ (0,x) ]
\in V(n,q)$. By definition,
$v = [(dv, \rho _{(n,q)}) v] $. So  $$  \left(
\begin{array}{c} dv \\ \rho_{(n,q)}v -x \end{array} \right) =
\left(
\begin{array}{cc} d               &  0  \\ \rho_{(n,q-1)}  & -d
\end{array} \right) \left( \begin{array}{c} t \\  y \end{array}
\right) = \left( \begin{array}{c} dt \\ \rho_{(n,q-1)}t - dy
\end{array} \right) $$ with $t \in N(n,q-1)^{^{n+1}}, y \in X^{^n}$. As
a consequence,
$$ ( \rho _{(n,q)} )^{^n}_{_*} [v-t] = [x] $$ and so $(
\rho_{(n,q)} )^{^n}_{_*}$ is an epimorphism. Lastly, put $
N(n+1,0)={\displaystyle
\lim_{\stackrel{\longrightarrow}{q}}} N(n,q),  $  $
\iota _{(n+1,0)} = {\displaystyle
\lim_{\stackrel{\longrightarrow}{q}}} \iota_{(n,q)} $,   $
\rho _{(n+1,0)} = {\displaystyle
\lim_{\stackrel{\longrightarrow}{q}}}
\rho _{(n,q)}   $ and let us verify the induction hypothesis.
Conditions  
$(i)_{(n+1,0)}$ and $(ii)_{(n+1,0)}$ are trivial. For
$(iii)_{(n+1,0)}$,  $( \rho _{(n+1,0)} )^{^i}_{_*} \colon
\hiru{H}{i}{N(n+1,0)} \longrightarrow \hiru{H}{i}{X}$ is  an
isomorphism for
$i = 0, \dots , n$ because all of the  $( \rho _{(n,q)} )^{^i}_{_*}$ are
isomorphisms by
$(iii)_{(n,q)}$ and $(
\rho _{(n+1,0)} )^{^{n+1}}_{_*}$ is a monomorphism, because if
$a
\in Z^{^{n+1}}N(n+1,0)$ is such that $[\rho_{(n+1,0)}a] = 0$, then, as
$a
\in N(n,q)$ for some $q$, we would have $\rho_{(n,q)}a = dx$ for
some $x
\in X$. So $v = [(a,x)] \in V(n,q+1)$ will kill the class
$[a]$ in $\hiru{H}{n+1}{N(n,q)}$ and therefore in
$\hiru{H}{n+1}{N(n+1,0)}$. \qed
\end{itemize}

\prg {\bf Proof of Theorem \ref{THII}.} We will confine ourselves to the case
$M=0$, since this is the only case we need in this paper. Let $N =A \otimes V$, and
assume we have already built
$$ 
\sigma_{(m,0)} \colon N(m,0) \to X
$$ 
in a such way that for all $m \leq n$

\medskip

$ 
\begin{array}{rl} 
(i)_m  & \rho\sigma_{(m,0)} = \Id_{N(m,0)},  \hbox{ and  } \\[.2cm] 
(ii)_m & \sigma_{(m,0)} \big|_{N(m',0)}  =
\sigma_{(m',0)} \hbox{ for all  } m'  \leq m.
\end{array} 
$

\bigskip

\nt We will define
$$ 
\sigma_{(n+1,0)} \colon N(n+1,0) \to X
$$
satisfying $(i)_{n+1}$, $(ii)_{n+1}$ and then the section of $\rho$ will be 
$$
\sigma = \lim_{\stackrel{\longrightarrow}{n}} \sigma_{(n,0)}.
$$
To do this, we will extend $\sigma_{(n,0)}$ to morphisms
$$ 
\sigma_{(n,q)} \colon N(n,q) \to X
$$
in such a way that for every $p\leq q$

\medskip

$ 
\begin{array}{rl} 
(i)_{(n,p)}  & \rho\sigma_{(n,p)} = \Id_{N(n,p)},  \hbox{ and  } \\[.2cm] 
(ii)_{(n,p)} & \sigma_{(n,p)}\big|_{N(n,p')}  =
\sigma_{(n,p')}, \hbox{ for all  } p'  \leq p.
\end{array} 
$

\bigskip

\nt Once we have these $\sigma_{(n,q)}$ for all $q\geq 0$, we will put
$$
\sigma_{(n+1,0)} = \lim_{\stackrel{\longrightarrow}{q}}  \sigma_{(n,q)},
$$
which will verify $(i)_{n+1}$ and $(ii)_{n+1}$.

So let us begin the construction of the $\sigma_{(n,q)}$. For $q=0$,
$\sigma_{(n,0)}$ exists by induction hypothesis. Assume we already have
$\sigma_{(n,p)}$ for all $p\leq q$. Then, consider the $A$-dg module
$$
X(n,q) = \Ima \left( \sigma_{(n,q)} \colon  N(n,q) \to X \right).
$$
Because $\rho \sigma_{(n,q)} = \Id_{N(n,q)}$, $\sigma_{(n,q)}$ is a monomorphism
and so $\sigma_{(n,q)} \colon N(n,q) \to X(n,q)$ is an isomorphism. If we apply the
Five Lemma to the long cohomology sequences of the following
commutative diagram of exact sequences of $A$-dg modules

\bigskip

$$ 
\begin{picture}(200,50)(0,0)

\put(-40,50){\makebox(0,0){$0$}}
\put(-40,0){\makebox(0,0){$0$}}

\put(20,50){\makebox(0,0){$X(n,q)$}}
\put(20,0){\makebox(0,0){$N(n,q)$}}

\put(80,50){\makebox(0,0){$X$}}
\put(80,0){\makebox(0,0){$N$}}

\put(140,50){\makebox(0,0){$X/X(n,q)$}}
\put(140,0){\makebox(0,0){$N/N(n,q)$}}

\put(200,50){\makebox(0,0){$0$}}
\put(200,0){\makebox(0,0){$0$}}

\put(20,42){\vector(0,-1){28}} 
\put(80,42){\vector(0,-1){28}} 
\put(140,42){\vector(0,-1){28}} 

\put(-30,0){\vector(1,0){25}} 
\put(-30,50){\vector(1,0){25}} 

\put(45,0){\vector(1,0){25}} 
\put(45,50){\vector(1,0){25}} 

\put(90,0){\vector(1,0){20}} 
\put(90,50){\vector(1,0){20}} 

\put(170,0){\vector(1,0){25}} 
\put(170,50){\vector(1,0){25}}

\put(35,28){\makebox(0,0){$\sigma^{-1}_{(n,q)}$}}
\put(87,28){\makebox(0,0){$\rho$}} 
\put(147,28){\makebox(0,0){$\bar\rho$}} 

\end{picture} 
$$ 

\bigskip

\nt we obtain that $\bar\rho$ is a {\em quis} of $A$-dg modules.  

Next, let $V(n,q+1)$ be an $\R$-vector space such that
$$
N(n,q+1) = N(n,q) \oplus \left( A \otimes V(n,q+1) \right)
$$
and let $j$ be the composition $ V \stackrel{i}{\hookrightarrow} N \to N/N(n,q)$.
Since $dV \subset N(n,q)$, we have $jV \subset Z(N/N(n,q))$, and so we have an
$\R$-linear morphism $H_j \colon V \to\hiru{H}{}{N/N(n,q)}$ which we can lift to
$\hiru{Z}{}{X/X(n,q)}$:

\bigskip

$$ 
\begin{picture}(200,50)(0,0)

\put(-20,0){\makebox(0,0){$V$}}
\put(80,0){\makebox(0,0){$\hiru{H}{}{N/N(n,q)}$}}
\put(200,60){\makebox(0,0){$\hiru{Z}{}{X/X(n,q)}$}}
\put(200,0){\makebox(0,0){$\hiru{H}{}{X/X(n,q)}$}}

\put(200,52){\vector(0,-1){38}} 
\put(200,52){\vector(0,-1){33}} 

\put(-10,0){\vector(1,0){50}}
\put(155,0){\vector(-1,0){35}} 

\put(-10,10){\vector(3,1){155}}

\put(20,8){\makebox(0,0){$H_j$}} 
\put(138,8){\makebox(0,0){$\cong$}} 
\put(77,50){\makebox(0,0){$\lambda$}}

\put(207,33){\makebox(0,0){$p$}} 

\end{picture} 
$$ 
\bigskip

\nt In fact, $N/N(n,q)$ is an $A$-dg module concentrated in degrees $\geq n$. So
$\hiru{B}{n}{N/N(n,q)} = 0$ and the previous diagram in degree $n$ is just

\bigskip

$$ 
\begin{picture}(200,50)(50,0)

\put(80,0){\makebox(0,0){$V$}}
\put(170,50){\makebox(0,0){$\hiru{Z}{n}{X/X(n,q)}$}}
\put(170,0){\makebox(0,0){$\hiru{Z}{n}{N/N(n,q)}$}}

\put(170,40){\vector(0,-1){26}}

\put(90,0){\vector(1,0){35}}

\put(90,10){\vector(1,1){35}}

\put(108,8){\makebox(0,0){$j$}} 
\put(108,45){\makebox(0,0){$\lambda$}}

\put(177,27){\makebox(0,0){$\bar\rho$}} 

\end{picture} 
$$ 
Next, consider the pull back

\bigskip

$$ 
\begin{picture}(200,60)(0, 0)

\put(40,0){\makebox(0,0){$X/X(n,q)$}}
\put(170,65){\makebox(0,0){$N$}}
\put(170,0){\makebox(0,0){$N/N(n,q)$}}
\put(40,65){\makebox(0,0){$\left( X/X(n,q) \right) \times_{_{ N/N(n,q)}} N$}}

\put(170,55){\vector(0,-1){46}} 
\put(80,0){\vector(1,0){55}}
\put(40,55){\vector(0,-1){46}} 
\put(105,65){\vector(1,0){50}}

\put(108,8){\makebox(0,0){$\bar\rho$}} 

\end{picture} 
$$ 

\bigskip

\nt and the morphism induced by $\pi \colon X \to X/X(n,q)$ and $\rho \colon X \to
N$,
$$
(\pi,\rho) \colon X \longrightarrow \left( X/X(n,q) \right) \times_{_{ N/N(n,q)}} N.
$$
It is an epimorphism: if $(\tilde{x},y) \in \left( X/X(n,q) \right) \times_{_{
N/N(n,q)}} N$ then $\rho x - y \in N(n,q)$. So, 
$$
\begin{array}{rcl}
(\pi,\rho ) (x - \sigma_{(n,q)} (\rho x-y)) &= &(\tilde{x} - 0, \rho x -(\rho x -y))
\\
&=& (\tilde{x},y).
\end{array}
$$
So, we can lift $(\lambda , i)$ to $X$

\bigskip

$$ 
\begin{picture}(200,50)(0,0)

\put(40,0){\makebox(0,0){$V$}}
\put(170,65){\makebox(0,0){$X$}}
\put(170,0){\makebox(0,0){$\left( X/X(n,q) \right) \times_{_{ N/N(n,q)}} N.$}}

\put(170,55){\vector(0,-1){46}} 
\put(170,55){\vector(0,-1){41}}

\put(50,0){\vector(1,0){55}}

\put(50,10){\vector(2,1){100}}

\put(78,8){\makebox(0,0){$(\lambda,i)$}} 
\put(108,50){\makebox(0,0){$f$}}

\put(187,37){\makebox(0,0){$(\pi,\rho)$}} 

\end{picture} 
$$ 

\bigskip

\nt Then $f$ and $\sigma_{(n,q)}$ induce the morphism of $A$-dg modules 
$
\sigma_{(n,q+1)} \colon N(n,q+1) \to X$ we were looking for:
$$
\sigma_{(n,q+1)} \big|_{N(n,q)} = \sigma_{(n,q)} \hbox { and } 
\sigma_{(n,q+1)} \big|_V =f.
$$
To see it is an $A$-dg morphism, it suffices, by definition of a Hirsch extension, to
verify that 
$$
\sigma_{(n,q)} d = df.
$$
So, let $v \in V$, then $\pi f v  = \lambda v \in \hiru{Z}{n}{X/X(n,q)}$. So, 
$0 = d \pi f v = \pi d fv$. Hence $dfv \in X(n,q) = \Ima \sigma_{(n,q)}$. Let
$\omega \in N(n,q)$ be such that $dfv = \sigma_{(n,q)} \omega$. Apply $\rho$ to
both sides of this equality and get, by the previous diagram,
$$
\rho d f v = d \rho fv = dv
$$
and, by induction hypothesis $(i)_{(n,q)}$,
$$
\rho \sigma_{(n,q)} \omega = \omega.
$$
So $dv = \omega$ , as we wanted.

Finally, $\sigma_{(n,q+1)}$ verifies $(i)_{(n,q+1)}$ and $(ii)_{(n,q+1)}$ by
construction. \qed

\prg {\bf Proof of Theorem \ref{abs}.}

This result generalizes to stratified spaces the result asserting that the
dgc algebra minimal model of a manifold can be computed using its deRham complex. The
proof is adapted from
\cite{Hal} using the same notation. It will be sufficient to prove
that the dgc algebra minimal model of
$\hiruv{\Om}{}{Z}$ is that of $A_{\R} (Z)$. For this purpose we construct a
commutative diagram

$$ \begin{picture}(300,60)(0,25)
\put(0,0){\makebox(0,0){$C(\underline{LS}(Z))$}}
\put(300,0){\makebox(0,0){$C(\underline{Sing}(Z))$}}

\put(30,60){\vector(1,0){70}} 
\put(260,0){\vector(-1,0){210}}
\put(265,60){\vector(-1,0){75}}

\put(0,45){\vector(0,-1){25}} 
\put(150,45){\vector(-3,-1){100}}
\put(300,45){\vector(0,-1){25}}

\put(310,32.5){\makebox(0,0){$\int_1$}}
\put(10,32.5){\makebox(0,0){$\int_2$}}
\put(140,32.5){\makebox(0,0){$\int_3$}}

\put(155,10){\makebox(0,0){$\rho_3$}}
\put(227.5,70){\makebox(0,0){$\rho_2$}}
\put(65.5,70){\makebox(0,0){$\rho_1$}}

\addtocounter{equation}{1}
\put(-90,32.5){\makebox(0,0){$(\theequation)$}}
\label{daig}
\put(0,60){\makebox(0,0){$\hiruv{\Om}{}{Z}$}}
\put(150,60){\makebox(0,0){$E(\underline{LS}(Z))$}}
\put(300,60){\makebox(0,0){$A_{\R} (Z)$}} 
\end{picture} 
$$

\vspace{1.5cm}

\nt where $\int_1, \int_2, \int_3, \hbox{ and } \rho_3$  induce cohomology
isomorphisms and
$\rho_2$ et $\rho_1$ are dgc algebra morphisms. This implies that $\rho_2$ and
$\rho_1$ are also dgc algebra quasi-isomorphisms and the Proposition is proved. We
construct the diagram in several steps.

\bigskip

\nt {\bf I. Unfolding of $\Delta$} (see \cite{S1}). The {\it unfolding} of the standard
simplex $\Delta$, relative to the decomposition
$\Delta = \Delta_0 * \cdots * \Delta_p$, is the map  
\\
$
\mu_\Delta\colon \wt{\Delta} = \bar{c}\Delta_0
\times
\cdots
\times \bar{c}\Delta_{p-1} \times  \Delta_p
\longrightarrow
\Delta $  defined by  $ \mu_\Delta ([x_0,t_0], \ldots, [x_{p-1},t_{p-1}],x_p) = 
t_0x_0+(1-t_0)t_1x_1+\cdots+(1-t_0)\cdots(1-t_{p-2})t_{p-1}x_{p-1}
+(1-t_0)\cdots(1-t_{p-1})x_p. $ Here
$\bar{c}\Delta_i$ denotes the closed cone  
$\Delta_i \times [0,1] \big/
\Delta_i \times \{ 0\}$ and $[x_i,t_i]$ a point of it.  This map is smooth and its
restriction $\mu_\Delta\colon
\inte(\wt{\Delta})  \longrightarrow \inte(\Delta) $ is a diffeomorphism (we write
$\inte(P) = P-\partial P$ the {\em interior} of the polyhedron
$P$). It sends a face $U$  of
$\wt{\Delta}$ on a face $V$ of $\Delta$ and the restriction
$\mu_\Delta \colon \inte(U) \to \inte(V)$ is a submersion.

This blow-up is compatible with face and degeneracy maps.  

{\em 1. Face}. Put 
$\delta_F \colon F \longrightarrow \Delta$ a codimension one face of
$\Delta$, the induced decomposition is 
$ F = \Delta_0*\cdots *\Delta_{j-1}* F_j *\Delta_{j+1} * \cdots*\Delta_p 
$ (we have written $\emptyset * X = X$).
 The {\it lifting} of $\delta_F$ is the map $\wt{\delta}_F \colon
\wt{F} \to \wt{\Delta}$ defined by: 

$
\wt{\delta}_F (z)= 
\left\{ 
\begin{array}{cl}  z & \hbox{if $F_j \neq \emptyset$}\\ (z_0, \ldots,
z_{j-1},\vartheta_j, z_{j+1} , \ldots , z_{p-1},x_{p}) & \hbox{if  $F_j =\emptyset$ ,
$j\neq p$}, \hbox{$\vartheta_j$ is the vertex of $\bar{c} \Delta_j$} \hbox{ and }\\  
& z = (z_0, \ldots, \widehat{z_j} , \ldots, z_{p-1}, x_{p})\\  
(z_0, \ldots, z_{p-2},[x_{p-1},1]) & \hbox{if $F_j=\emptyset$, $j=p$, and $z=(z_0,
\ldots, z_{p-2},x_{p-1})$}. 
\end{array}
\right.
$

 \nt This map is smooth (affine in barycentric coordinates), sends $\wt{F}$
isomorphically 
 on  a face of $\wt\Delta$ and verifies 
$\mu_\Delta \rond \wt{\delta}_F = \delta_F \rond
\mu_F$.  

\smallskip

{\em 2. Degeneracy}.   Put $\sigma \colon D = \Delta *
\{ P \}  \to
\Delta$ a degeneracy map with $\sigma_D ( P) = Q \in
\Delta_{j}$. The induced decomposition $D =
\Delta_0*\cdots *\Delta_{j-1}* (\Delta_{j} *\{ P\})*
\Delta_{j+1}*\cdots*\Delta_p)$.
 The {\it lifting} of
$\sigma$ is the map $\wt{\sigma}_D  \colon \wt{D} \to
\wt{\Delta}$ defined by:

$ 
\wt{\delta}_D (z)= \left\{ 
\begin{array}{cl} 
(z_0, \ldots ,  [tx_j + (1-t)Q,t_j]  ,\ldots ,z_{p-1},x_p) & \hbox{if
$j<p$ and $z$ is}\\ & (z_0, \ldots ,  [tx_j + (1-t)P,t_j]  ,\ldots ,z_{p-1},x_p), \\
 (z_0,\ldots,  z_{p-1},tx_p+ (1-t)Q) & \hbox{if $j=p$ and $z$ is } \\
& (z_0, \ldots ,
z_{p-1},tx_p+ (1-t)P).
\end{array}
\right.
$

  \nt This map is smooth (linear in barycentric coordinates) and verifies
$\mu_\Delta \rond \wt{\sigma}_D =
\sigma_D \rond \mu_D$.

On the boundary  $\partial\wt{\Delta}$ we find not only the blow-up
$\wt{\partial\Delta}$ of the boundary $
\partial\Delta$ of $\Delta$ but also  the faces  
$B_i =
\bar{c}\Delta_0 \times  \cdots \times
\bar{c}\Delta_{i-1} \times ( \Delta_i
\times \{ 1\}) \times \bar{c}\Delta_{i+1} \times 
\cdots \times
\bar{c}\Delta_{p-1} \times \Delta_p$ with  $i\in \{ 0,\ldots,p-2\}$ or $i=p-1$ and
$\dim \Delta_p >0$, which we shall call {\em bad faces}.  Notice that
$\dim \mu_\Delta ( B_i  )= \dim( \Delta_0 * \cdots *
\Delta_i)< n-1 =\dim B_i$.

\bigskip

\nt  {\bf II. The simplicial set\footnote{For the notions related with simplicial
sets, local systems, \ldots we refer the reader to
\cite{Hal}}}

$\underline{LS}(Z)$. On a stratified pseudomanifold it is not possible to
define smooth simplices directly  as in
\cite{Hal}. For this reason we introduce the notion of liftable simplices.Put  
$ Z=Z_{n =\dim R}
\supset \Sigma_Z = Z_{n-1}\supset\cdots
\supset Z_0 \supset Z_{-1} =\emptyset $ the filtration of $Z$, that is, $Z_i$ is the
union of strata with dimension smaller than $i$. A {\it liftable simplex} is a singular
simplex $\phii \colon
\Delta \to Z$ verifying the two following conditions. 
\begin{itemize} 
\item[{[LS1]}] {\em Each pull back $\phii^{-1}(Z_i)$ is a face of $\Delta$}.  
\end{itemize} 
 
Consider  $
\{i_0,\ldots,i_p\} = \{ i\in \{0,\ldots,n\} \ / \ 
\phii^{-1}(Z_i)\neq\phii^{-1}(Z_{i-1})\} $ and put
$\Delta_j$ the face of $\Delta$ with
$\phii^{-1}(Z_{i_j})= \phii^{-1}(Z_{i_j-1}) *
\Delta_j$. This defines on $\Delta$  the $\phii$-{\em decomposition}
$\Delta = \Delta_0 * \cdots * \Delta_p$.

\begin{itemize} 
\item[{[LS2]}] {\em There exists a smooth map
$\wt{\phii}
\colon \wt{\Delta} \to \wt{Z}$  with ${\cal L}_Z \rond
\wt{\phii }=
\phii \rond\mu_\Delta$, where the unfolding of $\Delta$ is taken relative to the
$\phii$-decomposition of $\Delta$}.
  \end{itemize}

The $\phii \rond\delta_F $-decomposition (resp.
$\phii\rond\sigma_D$-decomposition) is just  $F=\Delta_0*\cdots *\Delta_{j-1}*
F_j * \Delta_{j+1} *
\cdots*\Delta_p$ (resp. $D = \Delta_0*\cdots *\Delta_{j-1}* (\Delta_{j} *\{ P\})*
\Delta_{j+1}*\cdots*\Delta_p$). So, for a  liftable simplex 
$\phii$ the simplices $\partial_F(\phii) = \phii \rond \delta_F  $ and
$s_D (\phii) = \phii\rond \sigma_D  $ are  liftable simplices.  The face map
$\partial_F$ and the degeneracy map $s_D$ verify the usual compatibility
conditions. Put $\underline{LS}(Z)$ the family of liftable simplices. So, 
$( \underline{LS}(Z),\partial,s)$ define a simplicial set. 

 Notice that a liftable simplex
$\phii$ sends the interior of a face $A$ of $\Delta$ on a stratum $S$ of $Z$ and that
the restriction $\phii \colon int (A) \to S$ is smooth (using face maps we can
suppose
$A=\Delta$ and there we know that $\mu_\Delta \colon
\inte (\tilde\Delta) \to \inte\Delta$ is a diffeomorphism) . 

\bigskip

\nt  {\bf III. The local systems  $C$ and $E$}. The local system $C$ on
$\underline{Sing}(Z)$ (resp. $\underline{LS}(Z)$) is defined in
\cite[14.2]{Hal} in such a way that the space of global sections
$C(\underline{Sing}(Z))$ (resp. $C(\underline{LS}(Z))$) is the complex generated by
the  singular simplices (resp. liftable simplices) of $Z$.

Consider a simplex $\Delta$  endowed with the decomposition
$\Delta_0 *\cdots*\Delta_{p}$.  A {\em liftable form} on $\Delta$ is a family of
differential forms
$\eta = \{ \eta_A \in \hiru{\Om}{}{\inte(A)} \ / \ A \hbox{ face of }
\Delta
\}$ possessing  a common lifting $\wt{\eta}\in
\hiru{\Om}{}{\wt{\Delta}}$, that is: 
$$
\wt{\eta} =
\mu_\Delta^*\eta_{\mu_\Delta (H)} \hbox{   \ on  \ $\inte(H)$,   for each face $H$ of
$\wt{\Delta}$}.  
$$ The lifting
$\wt{\eta}$ is unique. Since $\mu_\Delta$ is an onto submersion with connected
fibers then the lifting forms are exactly the differential forms
$\om$ on
$\wt{\Delta}$
 verifying $\om (v,-)= d\om(v,-)=0$ for any vector of
$\wt{\Delta}$ with $(\mu_\Delta)_* (v)=0$, that is, the
 {\em basic forms} on $\wt\Delta$.  In the sequel we shall use both of points of
view.

A {\em liftable form} on $\{ \phii \colon
\Delta \to Z \} \in
\underline{LS}(Z)$ is a liftable form on $\Delta$ relative to its 
$\phii$-decom\-po\-si\-tion.
 For any $\phii\colon \Delta \to Z$  liftable simplex we shall write 
$E_{\phii}$ = \{liftable  forms on $\phii
$\}, which is a dgc algebra complex.

Consider $\delta_F \colon F \to \Delta$ a face map  and
$\sigma_D
\colon D \to \Delta$ a degeneracy map.   We define the face operator 
$\partial_F \colon E_{\phii}
\longrightarrow E_{\partial_F(\varphi)}$ and the degeneracy operator  
$ s_D \colon E_{\phii} \longrightarrow E_{s_D(\varphi)} $ by
 $\partial_F(\wt{\eta})=  \wt{\delta}_F^*\wt{\eta}$ and
 $s_D(\wt{\eta}) = \wt{\sigma}_D^*\wt{\eta}$ respectively.  These operators
verify the usual compatibility conditions and thus define a local system $E$ on
$\underline{LS}(Z)$. Notice that the space
$E(\underline{LS}(Z))$ of global sections of 
$E$ is a dgc algebra.

When $Z$ is a manifold endowed with the stratification $\{ Z\}$, liftable simplex
becomes smooth simplex and so $\underline{LS}(Z) = \underline{Sing}^\infty(Z)$.
Moreover, $E$ becomes the local system $A_\infty$ of $C^\infty$ differential
forms.

\bigskip

\nt  {\bf IV. Operators $\rho$ and $\int$.}

$\bullet$ The operator $\rho_3$ is just the inclusion, which makes sense since 
any liftable simplex is a singular one Proceeding as in
\cite{S1} one proves that this inclusion induces an isomorphism in homology and
therefore that  $\rho_3$ is a quasi-isomorphism in the category of graded
vector spaces.

\medskip

$\bullet$ The operator $\rho_2$ is defined as follows. For each
$\om  \in A_{\R}(Z)$ and  each liftable simplex $\phii$ of $Z$ we put 
$
\mu_\Delta^*\om_\varphi  
$  the lifting of 
$ (\rho_2 ( \om))_\varphi. 
$
 This operator is a dgc algebra operator.

\medskip

$\bullet$ The operator $\rho_1$ is defined as follows. For each
$\om\in\hiruv{\Om}{}{Z}$ and each liftable simplex $\phii$ of $Z$ we put  
$\wt{\phii}^*\wt{\om}$ the lifting of
$ (\rho_1(\om))_\varphi . $  This operator is a dgc algebra operator.

$\bullet$ The operator $\int_1$ is given by integration of differential forms on
simplices: 

$$ ( \int_1 \om )(\phii) = \int_\Delta
\om_\varphi   \hspace{.5cm} \hbox{  where  } 
\om \in A_{\R } (Z) \hbox{ and  } 
\phii
 \hbox{ is a singular simplex of } Z.
$$

\nt The deRham Theorem says that
$\int_1$ is a quasi-isomorphism in the category of graded vector spaces 
\cite{Hal}.

\medskip

$\bullet$ The operator $\int_2$ is given by integration of differential forms on
simplices:   

$$
 (\int_2 \omega)(\phii)  =
\int_{\inte(\Delta)} \phii^*\omega  =
\int_{\wt{\Delta}}
\wt{\phii}^*\wt{\omega} \hspace{.5cm}\hbox{  where  } \om
\in
\hiruv{\Om}{}{Z} \hbox{ and  } \phii \hbox{ is a liftable simplex of } Z.  
$$

\nt This operator is differential if we have 
${\displaystyle \int_{B_i}
\wt{\phii}^*\wt{\omega}}=0$. To prove this, we write $S$ the stratum of 
$Z$ containing $\phii (\inte(\mu_{\Delta}   ( B_i )))$; since
 $\wt{\phii}^*\wt{\omega} = \mu_{\Delta}^*
\phii^*
\om_S$  then we get ${\displaystyle \int_{B_i}
\wt{\phii}^*\wt{\omega}=0}$ because $\dim B_i > \dim \mu_\Delta (B_i)$.
Proceeding as in \cite{S1} one proves that $\int_2$ is a quasi-isomorphism in the
category of graded vector spaces.

\medskip

$\bullet$ The operator $\int_3$ is  is given by integration of differential forms on
simplices:  

$$ ( \int_3 \eta )(\phii)  = 
\int_{\inte(\Delta)} (\eta_\varphi)_\Delta = 
\int_{\wt{\Delta}}
\wt{\eta}_\varphi  \hspace{.5cm} \hbox{  where  } 
\eta \in E(\underline{LS}(Z)  ) \hbox{ and  }  
\phii
 \hbox{ is a liftable simplex of } Z.
$$

\nt This operator is differential since
${\displaystyle \int_{B_i} \wt{\eta}_\varphi=0}$; this comes from the equality 
$\wt{\eta}_\varphi = \mu_{\Delta}^* (\eta_\varphi)_{\mu_{\Delta} ( B_i )}$ on
$\inte(B_i)$. We already know that the local system $C$ on
$\underline{LS}(Z)$ is an extendable local system
\cite[Proposition 14.11]{Hal}; if we prove that $D$ is an extendable local system we
shall get  that $\int_3$ is a quasi-isomorphism in the category of graded vector
spaces
\cite[Theorem 12.27]{Hal}. This fact comes from the PoincarŽ Lemma and the
Extension Property we prove below.  For this purpose we fix a decomposition
$\Delta =
\Delta_0 * \cdots * \Delta_p$.

\bigskip

\nt {\bf PoincarŽ Lemma}.  {\em  $
\hiru{H}{}{ \{ \hbox{liftable forms on } \Delta \} } =
\R. $ }

\smallskip

\begin{itemize}
\item[]
 We prove $
\hiru{H}{}{ \{ \hbox{basic forms on } \wt{\Delta} \} } =
\R. $ Suppose first the case $\dim \Delta_p >0$. Put $\vartheta$ vertex of
$\Delta_p$. Consider the following maps:

\begin{itemize} 

\item[] $h_1 \colon \Delta \times [0,1] \to \Delta$ defined by
$h_1(r_0 x_0 +\cdots + r_p x_p,t) = r_0 x_0 +\cdots + r_{p-1} x_{p-1} +r_p
(t\vartheta +(1-t)x_p)$,

\item[] $h_2 \colon \wt{\Delta} \times [0,1] \to
\wt{\Delta}$ defined by $h_2(z_0, \ldots , z_{p-1},x_p,t) =
(z_0,\ldots,z_{p-1},t\vartheta +(1-t)x_p).$
\end{itemize} They are smooth homotopy maps between $\Delta$ (resp.
$\wt{\Delta}$) and $\Delta' = \Delta_0* \cdots *\Delta_{p-1}*\{ \vartheta \} $ (resp.
$\wt{\Delta'}$). Since $\mu_\Delta(h_2(z_0, \ldots , z_{p-1},x_p,t)) = 
h_1(\mu_\Delta(z_0, \ldots , z_{p-1},x_p),t)$ the basic forms are preserved and
$h_2$ induces a homotopy operator between $
 \{ \hbox{basic forms on } \wt{\Delta} \}  $ and $
 \{ \hbox{basic forms on } \wt{\Delta'} \}  . $

 Suppose $\dim \Delta_p=0$, that is,
$\Delta_p = \{
\vartheta \}$. Write $\nabla$ the simplex $\Delta$ endowed with the decomposition
$\Delta_0 *\cdots*\Delta_{p-2}*(\Delta_{p-1}*\{ \vartheta
\})$. Since the map $\tau \colon \wt{\nabla} \to 
\wt{\Delta}$, defined by  $\tau(z_0,
\ldots,z_{p-2},tx_{p-1}+(1-t)\vartheta) = (z_0,
\ldots,z_{p-2},[x_{p-1},t],\vartheta)$, is a diffeomorphism verifying $\mu_\Delta
\rond \tau = \mu_\nabla,
$ we get that the complexes  \newline $ \{ \hbox{basic forms on }
\wt{\Delta} \} $ and  $ \{
\hbox{basic forms on } \wt{\nabla} \} $ are isomorphic.

 Applying alternatively these two procedures we arrive to the case
$\Delta = \{
\vartheta \}$ and here it is clear that the cohomology of  $
\{ \hbox{basic forms on }
\wt{\Delta}\equiv \{ \vartheta \}\} $ is $\R$.
\end{itemize}

\bigskip

\nt {\bf Extension property}. {\em Each liftable form on 
$\partial \Delta$ possesses an extension to a liftable form on
$\Delta$}.

\smallskip

\begin{itemize}
\item[] Let $\eta =  
\{ \eta_A \in \hiru{\Om}{}{\inte(A)} 
\ / \ A \hbox{ face of } \partial\Delta \}$ be  a liftable form on $\partial \Delta$. 
Put  $\wt{\eta}\in \hiru{\Om}{}{\wt{\partial\Delta}}$ the lifting, which is  a basic
form verifying 
$
 \wt{\eta} =
\mu_\Delta^*\eta_{\mu_\Delta (H)} \hbox{ on  }\inte(H)
$  for each face $H$ of $\wt{\Delta}$.  Since the fibers of $\mu_\Delta \colon
\wt{\partial\Delta} \to \partial\Delta$ are not necessarily connected then we
cannot identify liftable forms with basic forms.

 Notice first that any  form on $\partial\wt{\Delta}$ possesses an extension to a
form on $\wt{\Delta}$. Moreover, if the form is basic then its extension is
necessarily basic (the restriction
$\mu_\Delta \colon \inte(\wt{\Delta}) \to \inte(\Delta)$ is a diffeomorphism). So,
it suffices to extend $\wt\eta$ to a basic form  defined on $\partial\wt{\Delta}$. 

Consider
$\delta_F \colon F \longrightarrow \Delta$ a face of
$\partial\Delta$. Using face maps one constructs a  smooth map 
$
\wt{\delta}_F \colon \wt{F} \to \wt{\Delta}
$ sending isomorphically  $\wt{F}$ on  a face of
$\wt{\partial\Delta}$ and verifying
$\mu_{\Delta} \rond \wt{\delta}_F = \delta_F \rond
\mu_F$.  Notice that in this case the map $\wt{\delta}_F$ is not necessarily unique
($\codim_\Delta F$ can be greater than 2). This equality
 implies that  $\eta|_F =  \{ \eta_A  \ / \ A \hbox{ face of } F \}$, restriction of
$\eta$ to $F$,  is a liftable form with lifting 
$\wt{\delta}_F^* \wt\eta$. Recall that this form is $\mu_F$-basic.

 For each $i \in \{ 0,\ldots,p-1\}$ we put $\nabla_i = \mu_\Delta (B_i)$, whose
induced decomposition is just  $\nabla_i  =
\Delta_0*\cdots*\Delta_i$.  Define the projection $pr_i \colon B_i \to
\wt{\nabla_i}$  by  
 $pr_i (z_0,\ldots, z_{i-1}, (x_i,1), z_{i+1}, \ldots,$ $z_{p-1} ,x_p)= (z_0,\ldots
z_{i-1},x_i)$. This map sends (the interior of) a face of $B_i$ on (the interior of) a
face of $\nabla_i$. We define on $B_i$ the form $\gamma_i = pr_i^*
{\wt\delta}_{\nabla_i}^*\wt{\eta}$, which is a basic form because $\mu_\Delta =
\mu_{\nabla_i} \rond pr_i$. Moreover,
$\gamma_i$ is the lifting of
$\eta |_{\nabla_i}$: if
$H$ is a face of
$B_i$ then 
$
\gamma_i =  pr_i^* \delta_{\nabla_i}^*\mu_{\Delta}^*\eta_{\mu_\Delta (H)} =
\mu^*_\Delta \eta_{\mu_\Delta (H)}$ on $\inte(H)$. Since the lifting is unique the
forms $\{ \gamma_0,\ldots,\gamma_{p-1}\}$  define a basic form on $\partial
\wt\Delta - \wt{\partial\Delta}$, the union of bad faces. Again, the uniqueness of
the lifting implies that this form coincides with
$\wt\eta$ on $ \wt{\partial\Delta}$. Therefore, the extension of
$\wt{\eta}$ is constructed.

\end{itemize}

\medskip

$\bullet$ {\em Commutativity}. One easily checks 
$\int_3 \rond \rho_1  = \int_2$ and
$\int_3 \rond \rho_2 = \rho_3 \rond \int_1$, that is,  diagram \refp{daig} is
commutative. Since $\int_1$, $\int_2$, $\int_3$, $\rho_3$ induce isomorphisms in
cohomology (as graded vector spaces) and $\rho_2$, $\rho_1$ are dgc algebra
morphisms then 
$\rho_2$, $\rho_1$ are dgc algebra quasi-isomorphisms. So, the dgc algebra minimal model
of $Z$ is that of $\hiruv{\Omega}{}{Z}$. \qed

\prg {\bf Proof of Theorem \ref{diagrama}.}

We shall say that a morphism $f \colon Z' \to Z$ between two stratified spaces
 is {\em good} if it satisfies 
  the two following conditions:
  
\begin{itemize} 
\item[{[P1]}] $f$ preserves controlled forms:  $f^* \omega \in  \hiruv{\Om}{}{Z'}$ for any $\omega \in  \hiruv{\Om}{}{Z}$
 \item[{[P2]}] $f$ preserves liftable simplices
 $f \rond \phii'$-decomposition = $\phii'$-decomposition  and 
\item[] $f \rond \phii'
\in
\underline{LS}(Z)$
 for any $\phii' \in\underline{LS}(Z')$.
\end{itemize}

Consider now the following diagram
$$ 
\begin{picture}(150,60)(0,0)
\put(-20,50){\makebox(0,0){$\hiruv{\Omega}{}{Z'}$}}
\put(50,50){\makebox(0,0){$E(\underline{LS}(Z'))$}}
\put(50,0){\makebox(0,0){$E(\underline{LS}(Z))$}}
\put(-20,0){\makebox(0,0){$\hiruv{\Omega}{}{Z}$}}
\put(120,50){\makebox(0,0){$A_{\R}(Z')$}}
\put(120,0){\makebox(0,0){$A_{\R}(Z)$}}

\put(-22,10){\vector(0,1){30}} 
\put(50,10){\vector(0,1){30}}
\put(120,10){\vector(0,1){30}}
\put(0,0){\vector(1,0){20}} 
\put(0,50){\vector(1,0){20}}
\put(100,0){\vector(-1,0){20}} 
\put(100,50){\vector(-1,0){20}}

\put(10,8){\makebox(0,0){$\rho_1$}}
\put(10,58){\makebox(0,0){$\rho'_1$}}
\put(90,8){\makebox(0,0){$\rho_2$}}
\put(90,58){\makebox(0,0){$\rho'_2$}}

\put(-12,25){\makebox(0,0){$f^*$}} 
\put(42,25){\makebox(0,0){$f^*$}}
\put(112,25){\makebox(0,0){$f^*$}} 

\end{picture} 
$$ where the pull backs are defined as follows.

For each $\omega \in A_{\R}(Z)$ and each liftable simple $\phii$ of $Z'$ we put 
$(f^*\om)_{\phii'} = 
\om _{f\rond \phii'}$. This operator is a dgc algebra morphism.
 
Since $f$ is good then $f^* \colon \hiruv{\Om}{}{Z} \to \hiruv{\Om}{}{Z'}$ is a
well-defined dgc algebra operator. For each $\eta \in E(\underline{LS}(Z))$ and
each liftable simple $\phii$ of $Z'$ we put 
$(f^*\eta)_{\phii} = 
\eta _{f\rond \phii}$, which makes sense since $f$ is good.
One easily checks $f^* \rond \rho'_1 = \rho_1 \rond f^*$ and $f^* \rond
\rho_2 = \rho'_2 \rond f^*$, which ends the proof
\qed

\prg {\bf Proof of Proposition \ref{good}.}

We proof first the following Lemma

\nt {\bf  Lemma}.  {\em 
Let $\phii \colon \Delta \to M$ be a simplex satisfying [LS1]. Then, the family of
strata of $M$ meeting $\Ima \phii$ is totally ordered.}

\begin{itemize}

\item[] We prove that if $F_1$, $F_2$ are two faces of $\Delta$ and $S_1,S_2$ two strata
of $M$ with $\phii (\inte (F_i))\cap S_i \not= \emptyset$, $i=1,2$, then 
$S_1\leq S_{2}$ or
$S_2\leq S_1$.

Since $\phii^{-1} (M_{\dim S_i})$ (resp. $\phii^{-1} (M_{\dim S_i-1})$) is a face of
$\Delta$ meeting $\inte (F_i)$ (resp. not containing $\inte (F_i)$) then $\inte (F_i)
\subset \phii^{-1}(M_{\dim S_i})$ (resp. $\inte (F_i) \cap 
\phii^{-1}(M_{\dim S_i-1}) =
\emptyset)$. So, $\inte (F_i) \subset \phii^{-1}(M_{\dim S_i}-M_{\dim S_i-1})$ and
by connectivity we get $\inte (F_i) \subset \phii^{-1} (S_i)$. Notice that this
implies
$\phii (F_i )\subset \overline{S_i}$.

Consider now $F_3$ the smallest face of $\Delta$ containing $F_1$ and $F_2$. Put
$S_3$ a stratum of $M$ with $\phii (\inte (F_3)) \cap S_3 \not= \emptyset$ (it
always exists!). From the previous paragraph we get $\phii( \inte 
(F_3)) \subset S_3$ and
$\phii (F_3) \subset \overline{S_3}$ and therefore $S_1 \cap \overline{S_3} \not=
\emptyset$, $S_2 \cap \overline{S_3} \not=\emptyset$ and so $S_1 \leq S_3$
 and $S_2 \leq S_3$.

Let us suppose $\dim S_1 \leq \dim S_2$. Since the face $\phii^{-1}(M_{\dim S_2})$
contains $\inte (F_1)$ and $\inte (F_2)$ then it also contains $F_3$ by minimality.
Thus
$S_3 \subset M_{\dim S_2}$ which gives $\dim S_3 \leq \dim S_2$ and therefore
$S_3=S_2$. Finally,  $S_1\leq S_2$. 
\end{itemize}

\nt Now, in ordre to verify Proposition \ref{good}, we need to verify conditions [P1] and [P2].
\smallskip

\nt $\bullet$ Projection $\pi$.

\smallskip

[P1] Since $\pi$ is a liftable morphism.

\smallskip

[P2] Let $\phii \colon \Delta \to M$ be a liftable simplex. Following the lemma the
family of strata meeting $\Ima \phii$ can be written 
$S_0 < S_1  < \cdots < S_{p-1} < S_p$.
 The $\phii$-decomposition of $\Delta$ is
$\Delta = \Delta_0 *
\cdots * \Delta_p$ with $\Delta_0 * \cdots *\Delta_i = \phii ^{-1}(S_0 \cup \cdots
\cup S_i)$. Since $\pi$ is a strict morphism then the family of strata of $B$
meeting $\Ima (\pi \rond \phii)$ is
$\pi(S_0)  < \pi(S_1 ) < \cdots < \pi(S_{p-1}) < \pi(S_p)$.
So, $\Delta_0 * \cdots *\Delta_i  = (\pi \rond \phii)^{-1}(\pi(S_0)  \cup 
\cdots \cup \pi(S_i))$ which implies that the $(\pi \rond \phii)$-decomposition of
$\Delta$ is the $\phii$-decomposition. The lifting of $\pi \rond \phii$ is just
$\wt{\pi} \rond \wt{\phii}$.

\smallskip

\nt $\bullet$ Inclusion $\iota$.

\smallskip

[P1] Since $\len (F) =0$ then we have $\hiruv{\Om}{}{F} = \hiru{\Om}{}{F}$. Notice
that $\iota^* \colon \hiruv{\Om}{}{B} \to \hiru{\Om}{}{F}$ is just $R_F$.

\smallskip

[P2] Consider $\phii \colon \Delta \to F$ a liftable simplex, that is, a smooth map
$\phii \colon \Delta \to S$ where $S$ is a fixed stratum. So, the
$\phii$-decomposition and the $(\iota \rond \phii)$-decomposition are just
$\Delta = \Delta$. On the other hand, since ${\cal L}_B \colon {\cal L}^{-1}_B (S) \to
S$ is a fibration and $\Delta$ is contractible then we can construct a smooth map
$\psi \colon \Delta \to \wt{Z}$ with $\mu_Z \rond \psi = \iota \rond \phii$. So, this
map is the lifting of $\iota \rond \phii$.
\qed


\begin{thebibliography}{99}

\bibitem{A} C. Allday: {\it On the rational homotopy of the fixed 
point 
set of  torus actions.} - Topology {\bf 17}(1978), 95-100.

\bibitem{AP} C. Allday and V. Puppe: {\sl Cohomological methods in 
transformation groups} - Cambridge Stud. Math. {\bf 32}, Cambridge 
Univ. Press, London and New York, 1993.

 \bibitem{BG} A.K. Bousfield and V.K.A.M. Gugenheim: {\it   On PL 
deRham theory
and rational homotopy type.} - Mem.  Amer. Math. Soc. {\bf 179}(1976).

\bibitem{Br1} G. Bredon: {\sl Sheaf theory\/}. - McGraw-Hill, New 
York, 1967.

\bibitem{Br2} G. Bredon: {\sl Introduction to compact  transformation 
groups\/}.- Pure and Appl. Math., Academic  Press, New York and London, 1972.

\bibitem{GKMP} R. Goresky, R. Kottwitz and MacPherson: {\it 
Equivariant cohomology, Koszul duality, and the localization theorem\/}. Invent. 
Math. {\bf 131}(1998), 25-83.

\bibitem{GHV} W. Greub, S. Halperin and R. Vanstone: {\sl 
Connections, curcture
and cohomology\/}. Pure and Apl. Math, Academic Press; New York and 
London, 1972.

\bibitem{GMo} P. A. Griffiths and J. W. Morgan: {\sl  Rational 
homotopy theory
and differential forms .} -  Progress in Mathematics {\bf 16}, 
BirkhŠuser,
Boston, Basel, Stuttgart, 1981.

\bibitem{GM} V.K.A.M. Gugenheim and J.P. May: {\it   On the theory and
applications of differential torsion products.} -  Mem.  Amer. Math. 
Soc. {\bf
142}(1974). 

\bibitem{Hal} S. Halperin: {\sl  Lectures on minimal models} -  Mem. 
Soc. Math. Fra.
 {\bf 9/10}(1983).

\bibitem{H-T} S. Halperin and D.  TanrŽ: {\it  Homotopie filtrŽe et 
fibrŽs
$C^{\infty}$\/}. - Illinois J. Maths. {\bf 34}(1990),  284-324.

\bibitem{HS} G. Hector and M. Saralegi: {\it Intersection homology of 
${\bf
S}^1$-actions.\/} -  Trans. Amer. Math. Soc. {\bf 338}(1993), 263-288.

\bibitem{KM} I. Kriz and J.P. May: {\it Operads, algebras and 
motives.\/} -  
AstŽrisque {\bf 233}(1995).

\bibitem{Mol} P. Molino: {\it Feuilletages de Lie ˆ feuilles 
denses.\/} -
SŽm. GŽo. Diff. 1982-83, Montpellier.

\bibitem{Na} V. Navarro Aznar: {\it Sur la connexion de Gauss-Manin 
en homotopie
rationnelle.\/} Ann. Ec. Nor. Sup. {\bf 26}(1993), 99-148. 

\bibitem{NR} M. Nicolau and A. Revent\'os: {\it On some geometrical 
properties
of Seifert bundles. \/} - Israel J. Math. {\bf 47}(1984), 323-334. 
 
\bibitem{Qui} D.G. Quillen: {\sl Homotopical algebra} - Springer 
Lect. Notes Math.  
{\bf  47}(1967).

\bibitem{R1} A. Roig: {\it Modles minimaux et foncteurs dŽrivŽs. \/} 
- J. Pure  Appl. Math. {\bf 91}(1994),  231-254. 
  
\bibitem{R2} A. Roig: {\it Formalizability of DG modules and 
morphisms of CDG
algebras. \/} - Illinois J. Math. {\bf 38}(1994), 434-451. 
      
\bibitem{R3} A. Roig: {\it Minimal resolutions and other minimal 
models.\/} - 
Publ. Mat. {\bf 37}(1993), 285-303. 
    
\bibitem{S1} M. Saralegi:  {\it Homological properties of stratified spaces.\/}
- Illinois J. Math. {\bf 38}(1994), 47-70. 

\bibitem{S3} M. Saralegi: {\it A Gysin sequence for semifree actions of 
$S^3$\/} - Proc. Amer. Math. Soc. {\bf 118}(1993), 1335-1345. 
   
\bibitem{S2} M. Saralegi: {\it Gysin sequences.\/} 
Analysis and geometry in foliated manifolds, Word Scientific Publishing (1995),
207-223.
                                                             
\bibitem{V} A. Verona : {\it Le ThŽorme de deRham  pour les
pr\'estratifications abstraites\/}. - C. R. Acad. Sci. Paris {\bf 273}(1971),
886-889.

 
\end{thebibliography}
\end{document}